\documentstyle[12pt,amssymb,epic,eepic,fullpage]{amsart}
\def\draft{n}
\theoremstyle{plain}

\newtheorem{theorem}{Theorem}
\newtheorem{fact}{Fact}
\newtheorem{proposition}{Proposition}[section]

\newtheorem{corollary}[proposition]{Corollary}

\theoremstyle{definition}
\newtheorem{definition}[proposition]{Definition}

\theoremstyle{remark}

\newtheorem{exercise}[proposition]{Exercise}

\newtheorem{remark}[proposition]{Remark}

\def\printname#1{
	\if\draft y
		\smash{\makebox[0pt]{\hspace{-0.5in}
			\raisebox{8pt}{\tt\tiny #1}}}
	\fi
}

\newcommand{\mathmode}[1]{$#1$}

\newlength{\standardunitlength}
\catcode`\@=11
\long\def\@makecaption#1#2{%
    \vskip 10pt
    \setbox\@tempboxa\hbox{
      \small\sf{\bfcaptionfont #1. }\ignorespaces #2}%
    \ifdim \wd\@tempboxa >\captionwidth {%
        \rightskip=\@captionmargin\leftskip=\@captionmargin
        \unhbox\@tempboxa\par}%
      \else
        \hbox to\hsize{\hfil\box\@tempboxa\hfil}%
    \fi}
\font\bfcaptionfont=cmssbx10 scaled \magstephalf
\newdimen\@captionmargin\@captionmargin=2\parindent
\newdimen\captionwidth\captionwidth=\hsize
\catcode`\@=12

\newcommand{\End}{\operatorname{End}}

\def\lbl#1{\label{#1}\printname{#1}}
\def\biblbl#1{\bibitem{#1}\printname{#1}}

\def\eqdef{\overset{\text{def}}{=}}

\newcommand{\G}[1]{{\cal G}_{#1}}

\newcommand{\tr}{\operatorname{tr}}

\newcommand{\calA}{{\cal A}}
\newcommand{\calD}{{\cal D}}
\newcommand{\calI}{{\cal I}}
\newcommand{\calK}{{\cal K}}
\newcommand{\calS}{{\cal S}}
\newcommand{\calW}{{\cal W}}
\newcommand{\EEPIC}[2]{
	\setlength{\unitlength}{#2\standardunitlength}
        \begin{array}{c}  \hspace{-1.7mm}
                \raisebox{-8pt}{#1}
                \hspace{-1.9mm}
        \end{array}
}

\def\Asl
{{
\begin{picture}(240,115)(0,-10)
\put(200.000,35.000){\arc{40.000}{3.1416}{6.2832}}
\path(0,0)(0,100)
\path(2.000,92.000)(0.000,100.000)(-2.000,92.000)
\path(80,0)(80,100)
\path(82.000,92.000)(80.000,100.000)(78.000,92.000)
\path(160,0)(160,100)
\path(162.000,92.000)(160.000,100.000)(158.000,92.000)
\path(240,0)(240,100)
\path(242.000,92.000)(240.000,100.000)(238.000,92.000)
\path(0,80)(80,80)
\path(80,20)(160,80)
\path(40,80)(160,20)
\path(130,35)(240,35)
\path(0,40)(80,40)
\end{picture}
}}
\def\Bridged
{{
\begin{picture}(705,106)(0,-10)
\path(0,30)	(2.915,26.504)
	(5.000,25.000)

\path(5,25)	(8.650,24.111)
	(12.921,23.984)
	(17.607,24.365)
	(22.500,25.000)
	(27.393,25.635)
	(32.079,26.016)
	(36.350,25.889)
	(40.000,25.000)

\path(40,25)	(43.161,23.257)
	(45.000,20.000)

\path(45,20)	(43.531,16.469)
	(40.000,15.000)

\path(40,15)	(36.469,16.469)
	(35.000,20.000)

\path(35,20)	(36.839,23.257)
	(40.000,25.000)

\path(40,25)	(43.650,25.889)
	(47.921,26.016)
	(52.607,25.635)
	(57.500,25.000)
	(62.393,24.365)
	(67.079,23.984)
	(71.350,24.111)
	(75.000,25.000)

\path(75,25)	(77.084,26.504)
	(80.000,30.000)

\put(55,30){\makebox(0,0)[lb]{
\raisebox{0pt}[0pt][0pt]{\shortstack[l]{\mathmode{\scriptstyle 2}}}}}
\put(20,30){\makebox(0,0)[lb]{
\raisebox{0pt}[0pt][0pt]{\shortstack[l]{\mathmode{\scriptstyle 1}}}}}
\path(140,30)	(142.916,26.504)
	(145.000,25.000)

\path(145,25)	(148.650,24.111)
	(152.921,23.984)
	(157.607,24.365)
	(162.500,25.000)
	(167.393,25.635)
	(172.079,26.016)
	(176.350,25.889)
	(180.000,25.000)

\path(180,25)	(183.161,23.257)
	(185.000,20.000)

\path(185,20)	(183.531,16.469)
	(180.000,15.000)

\path(180,15)	(176.469,16.469)
	(175.000,20.000)

\path(175,20)	(176.839,23.257)
	(180.000,25.000)

\path(180,25)	(183.650,25.889)
	(187.921,26.016)
	(192.607,25.635)
	(197.500,25.000)
	(202.393,24.365)
	(207.079,23.984)
	(211.350,24.111)
	(215.000,25.000)

\path(215,25)	(217.084,26.504)
	(220.000,30.000)

\put(195,30){\makebox(0,0)[lb]{
\raisebox{0pt}[0pt][0pt]{\shortstack[l]{\mathmode{\scriptstyle 2}}}}}
\put(160,30){\makebox(0,0)[lb]{
\raisebox{0pt}[0pt][0pt]{\shortstack[l]{\mathmode{\scriptstyle 1}}}}}
\path(220,30)	(222.916,26.504)
	(225.000,25.000)

\path(225,25)	(228.650,24.111)
	(232.921,23.984)
	(237.607,24.365)
	(242.500,25.000)
	(247.393,25.635)
	(252.079,26.016)
	(256.350,25.889)
	(260.000,25.000)

\path(260,25)	(263.161,23.257)
	(265.000,20.000)

\path(265,20)	(263.531,16.469)
	(260.000,15.000)

\path(260,15)	(256.469,16.469)
	(255.000,20.000)

\path(255,20)	(256.839,23.257)
	(260.000,25.000)

\path(260,25)	(263.650,25.889)
	(267.921,26.016)
	(272.607,25.635)
	(277.500,25.000)
	(282.393,24.365)
	(287.079,23.984)
	(291.350,24.111)
	(295.000,25.000)

\path(295,25)	(297.085,26.504)
	(300.000,30.000)

\put(275,30){\makebox(0,0)[lb]{
\raisebox{0pt}[0pt][0pt]{\shortstack[l]{\mathmode{\scriptstyle 2}}}}}
\put(240,30){\makebox(0,0)[lb]{
\raisebox{0pt}[0pt][0pt]{\shortstack[l]{\mathmode{\scriptstyle 1}}}}}
\path(300,30)	(302.915,26.504)
	(305.000,25.000)

\path(305,25)	(308.650,24.111)
	(312.921,23.984)
	(317.607,24.365)
	(322.500,25.000)
	(327.393,25.635)
	(332.079,26.016)
	(336.350,25.889)
	(340.000,25.000)

\path(340,25)	(343.161,23.257)
	(345.000,20.000)

\path(345,20)	(343.531,16.469)
	(340.000,15.000)

\path(340,15)	(336.469,16.469)
	(335.000,20.000)

\path(335,20)	(336.839,23.257)
	(340.000,25.000)

\path(340,25)	(343.650,25.889)
	(347.921,26.016)
	(352.607,25.635)
	(357.500,25.000)
	(362.393,24.365)
	(367.079,23.984)
	(371.350,24.111)
	(375.000,25.000)

\path(375,25)	(377.085,26.504)
	(380.000,30.000)

\put(355,30){\makebox(0,0)[lb]{
\raisebox{0pt}[0pt][0pt]{\shortstack[l]{\mathmode{\scriptstyle 2}}}}}
\put(320,30){\makebox(0,0)[lb]{
\raisebox{0pt}[0pt][0pt]{\shortstack[l]{\mathmode{\scriptstyle 1}}}}}
\put(130.000,-105.357){\arc{380.714}{4.0963}{5.3285}}
\put(240.000,44.643){\arc{80.714}{3.2747}{6.1501}}
\put(600.000,45.000){\arc{70.000}{3.1416}{6.2832}}
\put(600.000,45.000){\arc{90.000}{3.1416}{6.2832}}
\put(520.000,-12.778){\arc{205.556}{3.7386}{5.6862}}
\put(520.000,-17.857){\arc{195.714}{3.8391}{5.5856}}
\path(0,45)(80,45)
\path(72.000,43.000)(80.000,45.000)(72.000,47.000)
\path(140,45)(220,45)
\path(212.000,43.000)(220.000,45.000)(212.000,47.000)
\path(220,45)(300,45)
\path(292.000,43.000)(300.000,45.000)(292.000,47.000)
\path(300,45)(380,45)
\path(372.000,43.000)(380.000,45.000)(372.000,47.000)
\path(20,45)(20,50)
\path(60,45)(60,50)
\path(160,45)(160,50)
\path(200,45)(200,50)
\path(240,45)(240,50)
\path(280,45)(280,50)
\path(320,45)(320,50)
\path(360,45)(360,50)
\path(420,45)(435,45)
\path(445,45)(480,45)
\path(520,45)(555,45)
\path(565,45)(595,45)
\path(605,45)(635,45)
\path(645,45)(700,45)
\path(385,65)(415,65)
\path(407.000,63.000)(415.000,65.000)(407.000,67.000)
\path(420,45)	(416.377,42.294)
	(415.000,40.000)

\path(415,40)	(414.708,36.075)
	(415.483,31.861)
	(417.266,27.966)
	(420.000,25.000)

\path(420,25)	(423.157,22.955)
	(426.551,21.183)
	(430.148,19.683)
	(433.915,18.456)
	(437.818,17.502)
	(441.822,16.820)
	(445.894,16.411)
	(450.000,16.274)
	(454.106,16.411)
	(458.178,16.820)
	(462.182,17.502)
	(466.085,18.456)
	(469.852,19.683)
	(473.449,21.183)
	(476.843,22.955)
	(480.000,25.000)

\path(480,25)	(482.734,27.966)
	(484.517,31.861)
	(485.292,36.075)
	(485.000,40.000)

\path(485,40)	(483.623,42.294)
	(480.000,45.000)

\path(520,45)	(516.305,42.396)
	(515.000,40.000)

\path(515,40)	(515.182,37.166)
	(516.348,34.353)
	(520.000,30.000)

\path(520,30)	(524.134,27.146)
	(528.823,24.560)
	(531.352,23.364)
	(533.990,22.230)
	(536.730,21.157)
	(539.560,20.144)
	(542.472,19.190)
	(545.455,18.292)
	(548.501,17.449)
	(551.600,16.661)
	(554.742,15.925)
	(557.919,15.240)
	(561.119,14.604)
	(564.334,14.017)
	(567.554,13.477)
	(570.770,12.982)
	(573.972,12.530)
	(577.151,12.122)
	(580.297,11.754)
	(583.400,11.425)
	(586.451,11.135)
	(589.441,10.881)
	(592.360,10.663)
	(595.198,10.479)
	(597.946,10.326)
	(600.594,10.205)
	(605.554,10.049)
	(610.000,10.000)

\path(610,10)	(614.446,10.049)
	(619.406,10.205)
	(622.054,10.326)
	(624.802,10.479)
	(627.640,10.663)
	(630.559,10.881)
	(633.549,11.135)
	(636.600,11.425)
	(639.703,11.754)
	(642.849,12.122)
	(646.028,12.530)
	(649.230,12.982)
	(652.446,13.477)
	(655.666,14.017)
	(658.881,14.604)
	(662.081,15.240)
	(665.258,15.925)
	(668.400,16.661)
	(671.499,17.449)
	(674.545,18.292)
	(677.528,19.190)
	(680.440,20.144)
	(683.270,21.157)
	(686.010,22.230)
	(688.648,23.364)
	(691.177,24.560)
	(695.866,27.146)
	(700.000,30.000)

\path(700,30)	(703.652,34.353)
	(704.818,37.166)
	(705.000,40.000)

\path(705,40)	(703.695,42.396)
	(700.000,45.000)

\put(255,0){\makebox(0,0)[lb]{
\raisebox{0pt}[0pt][0pt]{\shortstack[l]{\mathmode{\scriptstyle S_3}}}}}
\put(35,0){\makebox(0,0)[lb]{
\raisebox{0pt}[0pt][0pt]{\shortstack[l]{\mathmode{\scriptstyle S_1}}}}}
\put(175,0){\makebox(0,0)[lb]{
\raisebox{0pt}[0pt][0pt]{\shortstack[l]{\mathmode{\scriptstyle S_1}}}}}
\put(335,0){\makebox(0,0)[lb]{
\raisebox{0pt}[0pt][0pt]{\shortstack[l]{\mathmode{\scriptstyle S_3}}}}}
\end{picture}
}}
\def\CDExample
{{
\begin{picture}(120,135)(0,-10)
\put(60,60){\ellipse{120}{120}}
\path(60,120)(120,60)
\path(100,100)(0,60)
\path(20,100)(60,0)
\path(20,20)(100,20)
\end{picture}
}}
\def\CombedBraid
{{
\begin{picture}(370,335)(0,-10)
\path(10,0)(10,30)
\path(10,40)(10,70)
\path(10,80)(10,160)
\path(10,170)(10,280)
\path(10,290)(10,320)
\path(12.000,312.000)(10.000,320.000)(8.000,312.000)
\path(50,120)(50,200)
\path(50,210)(50,230)
\path(50,240)(50,320)
\path(52.000,312.000)(50.000,320.000)(48.000,312.000)
\path(90,0)(90,95)
\path(90,225)(90,320)
\path(92.000,312.000)(90.000,320.000)(88.000,312.000)
\path(130,0)(130,215)
\path(85,272)(55,273)
\path(45,274)(15,275)
\path(55,290)(85,291)
\path(250,0)(250,320)
\path(252.000,312.000)(250.000,320.000)(248.000,312.000)
\path(290,0)(290,320)
\path(292.000,312.000)(290.000,320.000)(288.000,312.000)
\path(330,0)(330,320)
\path(332.000,312.000)(330.000,320.000)(328.000,312.000)
\path(370,0)(370,320)
\path(372.000,312.000)(370.000,320.000)(368.000,312.000)
\path(250,40)(290,40)
\path(250,80)(290,80)
\path(290,120)(330,120)
\path(250,160)(330,160)
\path(290,200)(330,200)
\path(290,240)(370,240)
\path(250,280)(370,280)
\path(50,0)	(50.000,3.266)
	(50.000,5.000)

\path(50,5)	(50.000,7.500)
	(50.000,10.000)

\path(50,10)	(50.000,12.500)
	(50.000,15.000)

\path(50,15)	(50.306,17.433)
	(50.000,20.000)

\path(50,20)	(46.019,23.798)
	(43.348,25.282)
	(40.451,26.539)
	(37.497,27.603)
	(34.657,28.510)
	(30.000,30.000)

\path(30,30)	(25.618,31.246)
	(22.886,31.873)
	(20.000,32.500)
	(17.114,33.127)
	(14.382,33.754)
	(10.000,35.000)

\path(10,35)	(7.496,35.660)
	(4.307,36.505)
	(0.000,40.000)

\path(0,40)	(0.889,42.356)
	(5.000,45.000)

\path(15,50)	(18.139,50.208)
	(21.058,50.421)
	(23.769,50.642)
	(26.281,50.872)
	(30.757,51.365)
	(34.573,51.917)
	(37.817,52.542)
	(40.578,53.254)
	(45.000,55.000)

\path(45,55)	(48.235,56.966)
	(50.000,60.000)

\path(50,60)	(47.070,64.305)
	(44.443,65.811)
	(41.385,66.989)
	(38.154,67.920)
	(35.007,68.679)
	(32.203,69.347)
	(30.000,70.000)

\path(30,70)	(25.618,71.246)
	(22.886,71.873)
	(20.000,72.500)
	(17.114,73.127)
	(14.382,73.754)
	(10.000,75.000)

\path(10,75)	(7.496,75.660)
	(4.307,76.505)
	(0.000,80.000)

\path(0,80)	(0.889,82.356)
	(5.000,85.000)

\path(15,90)	(18.139,90.208)
	(21.058,90.421)
	(23.769,90.642)
	(26.281,90.872)
	(30.757,91.365)
	(34.573,91.917)
	(37.817,92.542)
	(40.578,93.254)
	(45.000,95.000)

\path(45,95)	(47.825,97.148)
	(50.000,100.000)

\path(50,100)	(50.194,102.540)
	(50.000,105.000)

\path(50,105)	(50.000,106.734)
	(50.000,110.000)

\path(90,95)	(90.458,98.164)
	(90.000,100.000)

\path(90,100)	(86.019,103.798)
	(83.348,105.282)
	(80.451,106.539)
	(77.497,107.603)
	(74.657,108.510)
	(70.000,110.000)

\path(70,110)	(65.618,111.246)
	(62.886,111.873)
	(60.000,112.500)
	(57.114,113.127)
	(54.382,113.754)
	(50.000,115.000)

\path(50,115)	(47.496,115.660)
	(44.307,116.505)
	(40.000,120.000)

\path(40,120)	(40.889,122.356)
	(45.000,125.000)

\path(45,148)	(42.379,148.165)
	(39.941,148.321)
	(35.581,148.611)
	(31.847,148.874)
	(28.665,149.116)
	(23.667,149.564)
	(20.000,150.000)

\path(20,150)	(18.257,150.309)
	(15.000,151.000)

\path(5,155)	(0.889,157.644)
	(0.000,160.000)

\path(0,160)	(4.293,163.514)
	(7.480,164.362)
	(10.000,165.000)

\path(10,165)	(14.368,166.188)
	(17.080,166.780)
	(19.947,167.347)
	(22.821,167.870)
	(25.553,168.332)
	(30.000,169.000)

\path(30,169)	(32.201,169.242)
	(35.200,169.474)
	(39.350,169.719)
	(41.965,169.853)
	(45.000,170.000)

\path(55,170)	(58.438,170.050)
	(61.636,170.134)
	(64.607,170.253)
	(67.362,170.409)
	(72.273,170.842)
	(76.467,171.451)
	(80.039,172.254)
	(83.088,173.268)
	(88.000,176.000)

\path(88,176)	(90.000,180.000)

\path(90,180)	(88.000,184.000)

\path(88,184)	(84.176,186.252)
	(79.201,187.990)
	(76.610,188.672)
	(74.125,189.233)
	(70.000,190.000)

\path(70,190)	(67.771,190.238)
	(64.762,190.317)
	(60.622,190.238)
	(58.018,190.139)
	(55.000,190.000)

\path(45,195)	(40.889,197.644)
	(40.000,200.000)

\path(40,200)	(44.291,203.517)
	(47.478,204.365)
	(50.000,205.000)

\path(50,205)	(52.553,205.723)
	(55.701,206.411)
	(59.245,207.066)
	(62.988,207.694)
	(66.733,208.297)
	(70.281,208.880)
	(73.436,209.446)
	(76.000,210.000)

\path(76,210)	(78.914,210.496)
	(82.590,211.092)
	(86.222,212.142)
	(89.000,214.000)

\path(89,214)	(89.922,217.010)
	(90.000,220.000)

\path(90,220)	(90.053,221.737)
	(90.000,225.000)

\path(130,215)	(130.290,218.206)
	(130.000,220.000)

\path(130,220)	(127.825,222.852)
	(125.000,225.000)

\path(125,225)	(120.578,226.746)
	(117.817,227.458)
	(114.573,228.083)
	(110.757,228.635)
	(106.281,229.128)
	(103.769,229.358)
	(101.058,229.579)
	(98.139,229.792)
	(95.000,230.000)

\path(85,230)	(81.326,230.348)
	(77.911,230.682)
	(74.740,231.005)
	(71.802,231.319)
	(69.083,231.624)
	(66.571,231.924)
	(62.116,232.514)
	(58.334,233.101)
	(55.122,233.702)
	(52.378,234.330)
	(50.000,235.000)

\path(50,235)	(47.461,235.613)
	(44.276,236.463)
	(40.000,240.000)

\path(40,240)	(40.889,242.356)
	(45.000,245.000)

\path(55,250)	(58.844,250.594)
	(62.157,251.109)
	(64.995,251.551)
	(67.413,251.931)
	(71.212,252.538)
	(74.000,253.000)

\path(74,253)	(77.819,253.677)
	(80.861,254.235)
	(85.000,255.000)

\path(95,255)	(98.112,254.962)
	(101.007,254.949)
	(103.697,254.962)
	(106.193,255.003)
	(110.647,255.169)
	(114.456,255.457)
	(117.708,255.875)
	(120.490,256.432)
	(125.000,258.000)

\path(125,258)	(128.068,259.821)
	(130.000,263.000)

\path(130,263)	(128.609,266.067)
	(126.000,268.000)

\path(126,268)	(121.347,269.934)
	(118.474,270.638)
	(115.114,271.183)
	(111.178,271.580)
	(106.574,271.840)
	(103.993,271.923)
	(101.211,271.976)
	(98.218,272.002)
	(95.000,272.000)

\path(5,275)	(0.889,277.644)
	(0.000,280.000)

\path(0,280)	(4.267,283.549)
	(7.450,284.401)
	(10.000,285.000)

\path(10,285)	(12.390,285.632)
	(15.142,286.203)
	(18.359,286.724)
	(22.143,287.208)
	(26.596,287.665)
	(29.106,287.888)
	(31.822,288.109)
	(34.756,288.329)
	(37.922,288.550)
	(41.333,288.773)
	(45.000,289.000)

\path(96,291)	(98.831,291.249)
	(101.464,291.497)
	(106.173,292.001)
	(110.205,292.528)
	(113.641,293.090)
	(116.560,293.704)
	(119.040,294.384)
	(123.000,296.000)

\path(123,296)	(125.755,298.264)
	(128.000,301.000)

\path(128,301)	(129.213,303.939)
	(130.000,307.000)

\path(130,307)	(130.070,309.507)
	(130.000,312.000)

\path(130,312)	(130.000,314.775)
	(130.000,316.988)
	(130.000,320.000)

\path(132.000,312.000)(130.000,320.000)(128.000,312.000)
\path(55,130)	(58.218,129.998)
	(61.211,130.024)
	(63.993,130.077)
	(66.574,130.160)
	(71.178,130.420)
	(75.114,130.817)
	(78.474,131.362)
	(81.347,132.066)
	(86.000,134.000)

\path(86,134)	(88.548,135.992)
	(90.000,139.000)

\path(90,139)	(87.000,143.000)

\path(87,143)	(82.233,145.142)
	(79.275,145.954)
	(75.810,146.611)
	(71.744,147.130)
	(66.983,147.525)
	(64.312,147.680)
	(61.433,147.810)
	(58.333,147.916)
	(55.000,148.000)

\end{picture}
}}
\def\CommDiag
{{
\begin{picture}(209,109)(0,-10)
\path(75,80)(195,80)
\path(187.000,78.000)(195.000,80.000)(187.000,82.000)
\path(70,65)(120,20)
\path(112.716,23.865)(120.000,20.000)(115.392,26.838)
\path(200,65)(150,20)
\path(154.608,26.838)(150.000,20.000)(157.284,23.865)
\put(125,85){\makebox(0,0)[lb]{
\raisebox{0pt}[0pt][0pt]{\shortstack[l]{\mathmode{\ab}}}}}
\put(185,35){\makebox(0,0)[lb]{
\raisebox{0pt}[0pt][0pt]{\shortstack[l]{\mathmode{\bc}}}}}
\put(45,70){\makebox(0,0)[lb]{
\raisebox{0pt}[0pt][0pt]{\shortstack[l]{\mathmode{\a}}}}}
\put(200,70){\makebox(0,0)[lb]{
\raisebox{0pt}[0pt][0pt]{\shortstack[l]{\mathmode{\b}}}}}
\put(95,0){\makebox(0,0)[lb]{
\raisebox{0pt}[0pt][0pt]{\shortstack[l]{\mathmode{\c}}}}}
\put(0,35){\makebox(0,0)[lb]{
\raisebox{0pt}[0pt][0pt]{\shortstack[l]{\mathmode{\ac}}}}}
\end{picture}
}}

\def\CrossedIdentity
{{
\begin{picture}(40,55)(0,-10)
\path(0,0)(0,15)(15,15)
	(25,25)(40,25)(40,40)
\path(42.000,32.000)(40.000,40.000)(38.000,32.000)
\path(40,0)(40,15)(25,15)
	(15,25)(0,25)(0,40)
\path(2.000,32.000)(0.000,40.000)(-2.000,32.000)
\end{picture}
}}
\def\Crossedij
{{
\begin{picture}(180,77)(0,-10)
\put(165,20){\blacken\ellipse{2}{2}}
\put(165,20){\ellipse{2}{2}}
\put(135,20){\blacken\ellipse{2}{2}}
\put(135,20){\ellipse{2}{2}}
\put(150,20){\blacken\ellipse{2}{2}}
\put(150,20){\ellipse{2}{2}}
\put(105,20){\blacken\ellipse{2}{2}}
\put(105,20){\ellipse{2}{2}}
\put(75,20){\blacken\ellipse{2}{2}}
\put(75,20){\ellipse{2}{2}}
\put(90,20){\blacken\ellipse{2}{2}}
\put(90,20){\ellipse{2}{2}}
\put(45,20){\blacken\ellipse{2}{2}}
\put(45,20){\ellipse{2}{2}}
\put(30,20){\blacken\ellipse{2}{2}}
\put(30,20){\ellipse{2}{2}}
\put(15,20){\blacken\ellipse{2}{2}}
\put(15,20){\ellipse{2}{2}}
\path(180,0)(180,60)
\path(182.000,52.000)(180.000,60.000)(178.000,52.000)
\path(0,0)(0,60)
\path(2.000,52.000)(0.000,60.000)(-2.000,52.000)
\path(62.000,52.000)(60.000,60.000)(58.000,52.000)
\path(60,60)(60,45)(85,45)
	(95,35)(120,35)(120,0)
\path(122.000,52.000)(120.000,60.000)(118.000,52.000)
\path(120,60)(120,45)(95,45)
	(85,35)(60,35)(60,0)
\put(65,50){\makebox(0,0)[lb]{
\raisebox{0pt}[0pt][0pt]{\shortstack[l]{\mathmode{\scriptstyle i}}}}}
\put(125,50){\makebox(0,0)[lb]{
\raisebox{0pt}[0pt][0pt]{\shortstack[l]{\mathmode{\scriptstyle j}}}}}
\end{picture}
}}
\def\DeltaExample
{{
\begin{picture}(290,55)(0,-10)
\path(0,0)(0,40)
\path(2.000,32.000)(0.000,40.000)(-2.000,32.000)
\path(10,0)(10,40)
\path(12.000,32.000)(10.000,40.000)(8.000,32.000)
\path(40,0)(40,40)
\path(42.000,32.000)(40.000,40.000)(38.000,32.000)
\path(50,0)(50,40)
\path(52.000,32.000)(50.000,40.000)(48.000,32.000)
\path(80,0)(80,40)
\path(82.000,32.000)(80.000,40.000)(78.000,32.000)
\path(90,0)(90,40)
\path(92.000,32.000)(90.000,40.000)(88.000,32.000)
\path(120,0)(120,40)
\path(122.000,32.000)(120.000,40.000)(118.000,32.000)
\path(130,0)(130,40)
\path(132.000,32.000)(130.000,40.000)(128.000,32.000)
\path(160,0)(160,40)
\path(162.000,32.000)(160.000,40.000)(158.000,32.000)
\path(170,0)(170,40)
\path(172.000,32.000)(170.000,40.000)(168.000,32.000)
\path(200,0)(200,40)
\path(202.000,32.000)(200.000,40.000)(198.000,32.000)
\path(210,0)(210,40)
\path(212.000,32.000)(210.000,40.000)(208.000,32.000)
\path(240,0)(240,40)
\path(242.000,32.000)(240.000,40.000)(238.000,32.000)
\path(250,0)(250,40)
\path(252.000,32.000)(250.000,40.000)(248.000,32.000)
\path(280,0)(280,40)
\path(282.000,32.000)(280.000,40.000)(278.000,32.000)
\path(290,0)(290,40)
\path(292.000,32.000)(290.000,40.000)(288.000,32.000)
\path(0,20)(40,20)
\path(80,20)(130,20)
\path(170,20)(200,20)
\path(250,20)(290,20)
\put(215,15){\makebox(0,0)[lb]{
\raisebox{0pt}[0pt][0pt]{\shortstack[l]{\mathmode{+}}}}}
\put(135,15){\makebox(0,0)[lb]{
\raisebox{0pt}[0pt][0pt]{\shortstack[l]{\mathmode{+}}}}}
\put(55,15){\makebox(0,0)[lb]{
\raisebox{0pt}[0pt][0pt]{\shortstack[l]{\mathmode{+}}}}}
\end{picture}
}}
\def\LinearPath
{{
\begin{picture}(380,65)(0,-10)
\path(0,30)	(2.916,26.504)
	(5.000,25.000)

\path(5,25)	(8.650,24.111)
	(12.921,23.984)
	(17.607,24.365)
	(22.500,25.000)
	(27.393,25.635)
	(32.079,26.016)
	(36.350,25.889)
	(40.000,25.000)

\path(40,25)	(43.161,23.257)
	(45.000,20.000)

\path(45,20)	(43.531,16.469)
	(40.000,15.000)

\path(40,15)	(36.469,16.469)
	(35.000,20.000)

\path(35,20)	(36.839,23.257)
	(40.000,25.000)

\path(40,25)	(43.650,25.889)
	(47.921,26.016)
	(52.607,25.635)
	(57.500,25.000)
	(62.393,24.365)
	(67.079,23.984)
	(71.350,24.111)
	(75.000,25.000)

\path(75,25)	(77.084,26.504)
	(80.000,30.000)

\put(55,30){\makebox(0,0)[lb]{
\raisebox{0pt}[0pt][0pt]{\shortstack[l]{\mathmode{\scriptstyle 2}}}}}
\put(20,30){\makebox(0,0)[lb]{
\raisebox{0pt}[0pt][0pt]{\shortstack[l]{\mathmode{\scriptstyle 1}}}}}
\path(140,30)	(142.916,26.504)
	(145.000,25.000)

\path(145,25)	(148.650,24.111)
	(152.921,23.984)
	(157.607,24.365)
	(162.500,25.000)
	(167.393,25.635)
	(172.079,26.016)
	(176.350,25.889)
	(180.000,25.000)

\path(180,25)	(183.161,23.257)
	(185.000,20.000)

\path(185,20)	(183.531,16.469)
	(180.000,15.000)

\path(180,15)	(176.469,16.469)
	(175.000,20.000)

\path(175,20)	(176.839,23.257)
	(180.000,25.000)

\path(180,25)	(183.650,25.889)
	(187.921,26.016)
	(192.607,25.635)
	(197.500,25.000)
	(202.393,24.365)
	(207.079,23.984)
	(211.350,24.111)
	(215.000,25.000)

\path(215,25)	(217.085,26.504)
	(220.000,30.000)

\put(195,30){\makebox(0,0)[lb]{
\raisebox{0pt}[0pt][0pt]{\shortstack[l]{\mathmode{\scriptstyle 2}}}}}
\put(160,30){\makebox(0,0)[lb]{
\raisebox{0pt}[0pt][0pt]{\shortstack[l]{\mathmode{\scriptstyle 1}}}}}
\path(220,30)	(222.915,26.504)
	(225.000,25.000)

\path(225,25)	(228.650,24.111)
	(232.921,23.984)
	(237.607,24.365)
	(242.500,25.000)
	(247.393,25.635)
	(252.079,26.016)
	(256.350,25.889)
	(260.000,25.000)

\path(260,25)	(263.161,23.257)
	(265.000,20.000)

\path(265,20)	(263.531,16.469)
	(260.000,15.000)

\path(260,15)	(256.469,16.469)
	(255.000,20.000)

\path(255,20)	(256.839,23.257)
	(260.000,25.000)

\path(260,25)	(263.650,25.889)
	(267.921,26.016)
	(272.607,25.635)
	(277.500,25.000)
	(282.393,24.365)
	(287.079,23.984)
	(291.350,24.111)
	(295.000,25.000)

\path(295,25)	(297.085,26.504)
	(300.000,30.000)

\put(275,30){\makebox(0,0)[lb]{
\raisebox{0pt}[0pt][0pt]{\shortstack[l]{\mathmode{\scriptstyle 2}}}}}
\put(240,30){\makebox(0,0)[lb]{
\raisebox{0pt}[0pt][0pt]{\shortstack[l]{\mathmode{\scriptstyle 1}}}}}
\path(300,30)	(302.915,26.504)
	(305.000,25.000)

\path(305,25)	(308.650,24.111)
	(312.921,23.984)
	(317.607,24.365)
	(322.500,25.000)
	(327.393,25.635)
	(332.079,26.016)
	(336.350,25.889)
	(340.000,25.000)

\path(340,25)	(343.161,23.257)
	(345.000,20.000)

\path(345,20)	(343.531,16.469)
	(340.000,15.000)

\path(340,15)	(336.469,16.469)
	(335.000,20.000)

\path(335,20)	(336.839,23.257)
	(340.000,25.000)

\path(340,25)	(343.650,25.889)
	(347.921,26.016)
	(352.607,25.635)
	(357.500,25.000)
	(362.393,24.365)
	(367.079,23.984)
	(371.350,24.111)
	(375.000,25.000)

\path(375,25)	(377.085,26.504)
	(380.000,30.000)

\put(355,30){\makebox(0,0)[lb]{
\raisebox{0pt}[0pt][0pt]{\shortstack[l]{\mathmode{\scriptstyle 2}}}}}
\put(320,30){\makebox(0,0)[lb]{
\raisebox{0pt}[0pt][0pt]{\shortstack[l]{\mathmode{\scriptstyle 1}}}}}
\path(0,45)(80,45)
\path(72.000,43.000)(80.000,45.000)(72.000,47.000)
\path(140,45)(220,45)
\path(212.000,43.000)(220.000,45.000)(212.000,47.000)
\path(220,45)(300,45)
\path(292.000,43.000)(300.000,45.000)(292.000,47.000)
\path(300,45)(380,45)
\path(372.000,43.000)(380.000,45.000)(372.000,47.000)
\path(20,45)(20,50)
\path(60,45)(60,50)
\path(160,45)(160,50)
\path(200,45)(200,50)
\path(240,45)(240,50)
\path(280,45)(280,50)
\path(320,45)(320,50)
\path(360,45)(360,50)
\put(255,0){\makebox(0,0)[lb]{
\raisebox{0pt}[0pt][0pt]{\shortstack[l]{\mathmode{\scriptstyle S_3}}}}}
\put(35,0){\makebox(0,0)[lb]{
\raisebox{0pt}[0pt][0pt]{\shortstack[l]{\mathmode{\scriptstyle S_1}}}}}
\put(175,0){\makebox(0,0)[lb]{
\raisebox{0pt}[0pt][0pt]{\shortstack[l]{\mathmode{\scriptstyle S_1}}}}}
\put(335,0){\makebox(0,0)[lb]{
\raisebox{0pt}[0pt][0pt]{\shortstack[l]{\mathmode{\scriptstyle S_3}}}}}
\end{picture}
}}
\def\MaxExample
{{
\begin{picture}(432,113)(0,-10)
\put(40.000,72.000){\arc{40.000}{3.1416}{6.2832}}
\put(100.000,72.000){\arc{40.000}{3.1416}{6.2832}}
\put(55.000,12.000){\arc{60.000}{3.1416}{6.2832}}
\put(85.000,12.000){\arc{60.000}{3.1416}{6.2832}}
\put(240.000,72.000){\arc{50.000}{3.1416}{6.2832}}
\put(240.000,72.000){\arc{30.000}{3.1416}{6.2832}}
\put(300.000,72.000){\arc{50.000}{3.1416}{6.2832}}
\put(300.000,72.000){\arc{30.000}{3.1416}{6.2832}}
\put(255.000,12.000){\arc{50.000}{3.1416}{6.2832}}
\put(285.000,12.000){\arc{50.000}{3.1416}{6.2832}}
\put(255.000,12.000){\arc{70.000}{3.1416}{6.2832}}
\put(285.000,12.000){\arc{70.000}{3.1416}{6.2832}}
\path(0,72)(140,72)
\path(132.000,70.000)(140.000,72.000)(132.000,74.000)
\path(0,12)(140,12)
\path(132.000,10.000)(140.000,12.000)(132.000,14.000)
\path(150,82)(190,82)
\path(182.000,80.000)(190.000,82.000)(182.000,84.000)
\path(150,22)(190,22)
\path(182.000,20.000)(190.000,22.000)(182.000,24.000)
\path(225,72)(255,72)
\path(265,72)(275,72)
\path(285,72)(315,72)
\path(230,12)(250,12)
\path(260,12)(280,12)
\path(290,12)(310,12)
\path(215,72)	(211.325,72.435)
	(208.589,72.581)
	(205.000,72.000)

\path(205,72)	(201.839,70.257)
	(200.000,67.000)

\path(200,67)	(201.839,63.743)
	(205.000,62.000)

\path(205,62)	(207.629,61.259)
	(210.733,60.912)
	(214.146,60.866)
	(217.700,61.032)
	(221.230,61.320)
	(224.567,61.637)
	(227.546,61.894)
	(230.000,62.000)

\path(230,62)	(233.961,62.000)
	(238.813,62.000)
	(241.480,62.000)
	(244.258,62.000)
	(247.110,62.000)
	(250.000,62.000)
	(252.890,62.000)
	(255.742,62.000)
	(258.520,62.000)
	(261.187,62.000)
	(266.039,62.000)
	(270.000,62.000)

\path(270,62)	(273.466,62.000)
	(277.711,62.000)
	(282.476,62.000)
	(284.972,62.000)
	(287.500,62.000)
	(290.028,62.000)
	(292.524,62.000)
	(297.289,62.000)
	(301.534,62.000)
	(305.000,62.000)

\path(305,62)	(307.945,61.873)
	(311.520,61.565)
	(315.525,61.184)
	(319.760,60.839)
	(324.025,60.639)
	(328.120,60.694)
	(331.845,61.111)
	(335.000,62.000)

\path(335,62)	(338.161,63.743)
	(340.000,67.000)

\path(340,67)	(338.161,70.257)
	(335.000,72.000)

\path(335,72)	(331.411,72.581)
	(328.675,72.435)
	(325.000,72.000)

\path(220,12)	(216.269,12.235)
	(213.515,12.314)
	(210.000,12.000)

\path(210,12)	(207.405,11.455)
	(204.226,10.597)
	(200.000,7.000)

\path(200,7)	(201.703,3.932)
	(205.000,2.000)

\path(205,2)	(207.629,1.259)
	(210.733,0.912)
	(214.146,0.866)
	(217.700,1.032)
	(221.230,1.320)
	(224.567,1.637)
	(227.546,1.894)
	(230.000,2.000)

\path(230,2)	(233.961,2.000)
	(238.813,2.000)
	(241.480,2.000)
	(244.258,2.000)
	(247.110,2.000)
	(250.000,2.000)
	(252.890,2.000)
	(255.742,2.000)
	(258.520,2.000)
	(261.187,2.000)
	(266.039,2.000)
	(270.000,2.000)

\path(270,2)	(273.466,2.000)
	(277.711,2.000)
	(282.476,2.000)
	(284.972,2.000)
	(287.500,2.000)
	(290.028,2.000)
	(292.524,2.000)
	(297.289,2.000)
	(301.534,2.000)
	(305.000,2.000)

\path(305,2)	(307.945,1.873)
	(311.520,1.565)
	(315.525,1.184)
	(319.760,0.839)
	(324.025,0.639)
	(328.120,0.694)
	(331.845,1.111)
	(335.000,2.000)

\path(335,2)	(338.297,3.932)
	(340.000,7.000)

\path(340,7)	(335.774,10.597)
	(332.595,11.455)
	(330.000,12.000)

\path(330,12)	(326.485,12.314)
	(323.731,12.235)
	(320.000,12.000)

\put(355,17){\makebox(0,0)[lb]{ \raisebox{0pt}[0pt][0pt]{\shortstack[l]{{\rm
(not maximal)}}}}}
\put(355,77){\makebox(0,0)[lb]{ \raisebox{0pt}[0pt][0pt]{\shortstack[l]{{\rm
(maximal)}}}}}
\end{picture}
}}
\def\Omegaij
{{
\begin{picture}(180,77)(0,-10)
\put(165,20){\blacken\ellipse{2}{2}}
\put(165,20){\ellipse{2}{2}}
\put(135,20){\blacken\ellipse{2}{2}}
\put(135,20){\ellipse{2}{2}}
\put(150,20){\blacken\ellipse{2}{2}}
\put(150,20){\ellipse{2}{2}}
\put(105,20){\blacken\ellipse{2}{2}}
\put(105,20){\ellipse{2}{2}}
\put(75,20){\blacken\ellipse{2}{2}}
\put(75,20){\ellipse{2}{2}}
\put(90,20){\blacken\ellipse{2}{2}}
\put(90,20){\ellipse{2}{2}}
\put(45,20){\blacken\ellipse{2}{2}}
\put(45,20){\ellipse{2}{2}}
\put(30,20){\blacken\ellipse{2}{2}}
\put(30,20){\ellipse{2}{2}}
\put(15,20){\blacken\ellipse{2}{2}}
\put(15,20){\ellipse{2}{2}}
\path(180,0)(180,60)
\path(182.000,52.000)(180.000,60.000)(178.000,52.000)
\path(120,0)(120,60)
\path(122.000,52.000)(120.000,60.000)(118.000,52.000)
\path(60,0)(60,60)
\path(62.000,52.000)(60.000,60.000)(58.000,52.000)
\path(0,0)(0,60)
\path(2.000,52.000)(0.000,60.000)(-2.000,52.000)
\path(60,40)(120,40)
\put(65,50){\makebox(0,0)[lb]{
\raisebox{0pt}[0pt][0pt]{\shortstack[l]{\mathmode{\scriptstyle i}}}}}
\put(125,50){\makebox(0,0)[lb]{
\raisebox{0pt}[0pt][0pt]{\shortstack[l]{\mathmode{\scriptstyle j}}}}}
\end{picture}
}}
\def\PBExample
{{
\begin{picture}(280,295)(0,-10)
\put(20,80){\blacken\ellipse{6}{6}}
\put(20,80){\ellipse{6}{6}}
\put(20,160){\blacken\ellipse{6}{6}}
\put(20,160){\ellipse{6}{6}}
\put(60,200){\blacken\ellipse{6}{6}}
\put(60,200){\ellipse{6}{6}}
\put(20,240){\blacken\ellipse{6}{6}}
\put(20,240){\ellipse{6}{6}}
\path(200,0)(200,280)
\path(200,0)(200,280)
\path(202.000,272.000)(200.000,280.000)(198.000,272.000)
\path(240,0)(240,280)
\path(240,0)(240,280)
\path(242.000,272.000)(240.000,280.000)(238.000,272.000)
\path(280,0)(280,280)
\path(280,0)(280,280)
\path(282.000,272.000)(280.000,280.000)(278.000,272.000)
\path(120,140)(160,140)
\path(120,140)(160,140)
\path(152.000,138.000)(160.000,140.000)(152.000,142.000)
\path(200,80)(280,80)
\path(240,160)(280,160)
\path(200,200)(280,200)
\path(200,240)(240,240)
\path(0,0)	(0.000,4.045)
	(0.000,7.531)
	(0.000,10.517)
	(0.000,13.062)
	(0.000,17.062)
	(0.000,20.000)

\path(0,20)	(-0.080,24.400)
	(-0.150,27.116)
	(-0.213,29.982)
	(-0.249,32.850)
	(-0.239,35.573)
	(0.000,40.000)

\path(0,40)	(0.730,44.502)
	(1.281,47.242)
	(1.926,50.114)
	(2.643,52.971)
	(3.409,55.668)
	(5.000,60.000)

\path(5,60)	(7.898,64.789)
	(9.883,67.554)
	(12.054,70.397)
	(14.281,73.190)
	(16.433,75.806)
	(20.000,80.000)

\path(20,80)	(24.276,84.550)
	(26.993,87.279)
	(29.884,90.128)
	(32.790,92.961)
	(35.550,95.642)
	(40.000,100.000)

\path(40,100)	(42.203,102.203)
	(45.204,105.203)
	(49.352,109.352)
	(51.966,111.966)
	(55.000,115.000)

\path(80,0)	(80.159,3.748)
	(80.212,6.508)
	(80.000,10.000)

\path(80,10)	(79.251,13.432)
	(78.033,17.651)
	(76.548,21.796)
	(75.000,25.000)

\path(75,25)	(72.011,29.109)
	(67.990,33.655)
	(63.255,38.459)
	(60.719,40.901)
	(58.123,43.339)
	(55.508,45.751)
	(52.914,48.114)
	(47.945,52.604)
	(43.534,56.626)
	(40.000,60.000)

\path(40,60)	(36.448,63.343)
	(31.993,67.289)
	(29.532,69.441)
	(26.970,71.685)
	(24.350,74.003)
	(21.712,76.375)
	(19.099,78.783)
	(16.552,81.206)
	(11.825,86.024)
	(7.863,90.673)
	(5.000,95.000)

\path(5,95)	(3.917,97.447)
	(2.943,100.435)
	(2.089,103.790)
	(1.366,107.339)
	(0.785,110.908)
	(0.356,114.323)
	(0.091,117.412)
	(0.000,120.000)

\path(0,120)	(0.091,122.588)
	(0.356,125.677)
	(0.785,129.092)
	(1.366,132.661)
	(2.089,136.210)
	(2.943,139.565)
	(3.917,142.553)
	(5.000,145.000)

\path(5,145)	(7.873,149.316)
	(11.861,153.938)
	(16.621,158.722)
	(19.183,161.129)
	(21.809,163.523)
	(24.457,165.884)
	(27.084,168.195)
	(29.647,170.439)
	(32.103,172.596)
	(36.522,176.579)
	(40.000,180.000)

\path(40,180)	(43.636,183.791)
	(48.121,188.375)
	(50.582,190.894)
	(53.135,193.527)
	(55.741,196.247)
	(58.360,199.025)
	(60.952,201.834)
	(63.477,204.645)
	(65.896,207.431)
	(68.169,210.163)
	(70.256,212.814)
	(72.116,215.356)
	(75.000,220.000)

\path(75,220)	(76.679,224.332)
	(77.449,227.029)
	(78.152,229.886)
	(78.774,232.757)
	(79.299,235.497)
	(80.000,240.000)

\path(80,240)	(80.239,244.427)
	(80.249,247.150)
	(80.213,250.018)
	(80.150,252.884)
	(80.080,255.600)
	(80.000,260.000)

\path(80,260)	(80.000,262.938)
	(80.000,266.938)
	(80.000,269.483)
	(80.000,272.469)
	(80.000,275.955)
	(80.000,280.000)

\path(82.000,272.000)(80.000,280.000)(78.000,272.000)
\path(65,45)	(67.160,47.917)
	(69.001,50.449)
	(71.845,54.531)
	(73.767,57.597)
	(75.000,60.000)

\path(75,60)	(76.594,64.302)
	(77.415,66.972)
	(78.193,69.803)
	(78.885,72.660)
	(79.446,75.403)
	(80.000,80.000)

\path(80,80)	(79.951,82.608)
	(79.706,85.705)
	(79.279,89.122)
	(78.689,92.686)
	(77.949,96.228)
	(77.078,99.574)
	(76.089,102.555)
	(75.000,105.000)

\path(75,105)	(72.075,108.847)
	(67.838,112.982)
	(63.431,116.876)
	(60.000,120.000)

\path(60,120)	(55.594,124.406)
	(52.871,127.129)
	(50.000,130.000)
	(47.129,132.871)
	(44.406,135.594)
	(40.000,140.000)

\path(40,140)	(36.448,143.343)
	(31.993,147.289)
	(29.532,149.441)
	(26.970,151.685)
	(24.350,154.003)
	(21.712,156.375)
	(19.099,158.783)
	(16.552,161.206)
	(11.825,166.024)
	(7.863,170.673)
	(5.000,175.000)

\path(5,175)	(3.911,177.445)
	(2.922,180.426)
	(2.051,183.772)
	(1.311,187.314)
	(0.721,190.878)
	(0.294,194.295)
	(0.049,197.392)
	(0.000,200.000)

\path(0,200)	(0.508,204.597)
	(1.029,207.341)
	(1.684,210.198)
	(2.441,213.031)
	(3.267,215.700)
	(5.000,220.000)

\path(5,220)	(7.520,223.990)
	(11.020,228.407)
	(15.188,233.091)
	(19.706,237.885)
	(24.262,242.631)
	(28.539,247.170)
	(32.223,251.346)
	(35.000,255.000)

\path(35,255)	(37.804,259.820)
	(39.158,262.628)
	(40.000,265.000)

\path(40,265)	(40.353,267.255)
	(40.470,270.272)
	(40.353,274.404)
	(40.206,276.997)
	(40.000,280.000)

\path(42.577,272.167)(40.000,280.000)(38.587,271.876)
\path(65,125)	(67.160,127.917)
	(69.001,130.449)
	(71.845,134.531)
	(73.767,137.597)
	(75.000,140.000)

\path(75,140)	(76.594,144.302)
	(77.415,146.972)
	(78.193,149.803)
	(78.885,152.660)
	(79.446,155.403)
	(80.000,160.000)

\path(80,160)	(79.951,162.608)
	(79.706,165.705)
	(79.279,169.122)
	(78.689,172.686)
	(77.949,176.228)
	(77.078,179.574)
	(76.089,182.555)
	(75.000,185.000)

\path(75,185)	(72.137,189.327)
	(68.175,193.976)
	(63.448,198.794)
	(60.901,201.217)
	(58.288,203.625)
	(55.650,205.997)
	(53.030,208.315)
	(50.468,210.559)
	(48.007,212.711)
	(43.552,216.657)
	(40.000,220.000)

\path(40,220)	(36.466,223.374)
	(32.055,227.396)
	(27.086,231.886)
	(24.492,234.249)
	(21.877,236.661)
	(19.281,239.099)
	(16.745,241.541)
	(12.010,246.345)
	(7.989,250.891)
	(5.000,255.000)

\path(5,255)	(3.452,258.204)
	(1.967,262.349)
	(0.749,266.568)
	(0.000,270.000)

\path(0,270)	(-0.211,273.492)
	(-0.159,276.252)
	(0.000,280.000)

\path(1.606,271.912)(0.000,280.000)(-2.389,272.107)
\path(55,35)	(52.927,32.004)
	(51.147,29.416)
	(48.352,25.289)
	(46.381,22.268)
	(45.000,20.000)

\path(45,20)	(42.252,15.160)
	(40.864,12.365)
	(40.000,10.000)

\path(40,10)	(39.686,6.485)
	(39.765,3.731)
	(40.000,0.000)

\end{picture}
}}
\def\PairPath
{{
\begin{picture}(645,253)(0,-10)
\dottedline{5}(30,42)(80,62)
\path(73.315,57.172)(80.000,62.000)(71.829,60.886)
\thicklines
\path(135,147)(140,142)(140,62)(135,57)
\path(165,147)(160,142)(160,62)(165,57)
\path(215,147)(220,142)(220,62)(215,57)
\path(245,147)(240,142)(240,62)(245,57)
\path(295,147)(300,142)(300,62)(295,57)
\path(325,147)(320,142)(320,62)(325,57)
\thinlines
\dottedline{5}(310,203)(270,188)
\path(276.788,192.682)(270.000,188.000)(278.193,188.936)
\path(467.000,144.000)(465.000,152.000)(463.000,144.000)
\path(465,152)(465,52)
\path(477.000,144.000)(475.000,152.000)(473.000,144.000)
\path(475,152)(475,52)
\path(627.000,144.000)(625.000,152.000)(623.000,144.000)
\path(625,152)(625,52)
\path(637.000,144.000)(635.000,152.000)(633.000,144.000)
\path(635,152)(635,52)
\path(465,152)(465,192)
\thicklines
\path(455,147)(460,142)(460,62)(455,57)
\path(485,147)(480,142)(480,62)(485,57)
\path(535,147)(540,142)(540,62)(535,57)
\path(565,147)(560,142)(560,62)(565,57)
\path(615,147)(620,142)(620,62)(615,57)
\path(645,147)(640,142)(640,62)(645,57)
\thinlines
\path(358.000,104.000)(350.000,102.000)(358.000,100.000)
\path(350,102)(390,102)
\path(382.000,100.000)(390.000,102.000)(382.000,104.000)
\path(115,32)	(117.793,31.581)
	(120.537,31.171)
	(123.232,30.771)
	(125.880,30.380)
	(128.480,29.998)
	(131.034,29.626)
	(136.004,28.909)
	(140.796,28.229)
	(145.412,27.585)
	(149.860,26.976)
	(154.142,26.403)
	(158.265,25.864)
	(162.232,25.360)
	(166.050,24.889)
	(169.721,24.452)
	(173.252,24.047)
	(176.647,23.675)
	(179.911,23.335)
	(183.049,23.026)
	(186.065,22.749)
	(188.965,22.502)
	(191.753,22.285)
	(194.433,22.098)
	(197.012,21.940)
	(199.492,21.811)
	(204.181,21.637)
	(208.538,21.573)
	(212.601,21.614)
	(216.409,21.758)
	(220.000,22.000)

\path(220,22)	(223.674,22.326)
	(227.793,22.740)
	(232.289,23.248)
	(237.098,23.860)
	(239.598,24.207)
	(242.152,24.582)
	(244.750,24.987)
	(247.385,25.422)
	(250.049,25.889)
	(252.732,26.388)
	(255.427,26.921)
	(258.125,27.488)
	(260.818,28.090)
	(263.498,28.728)
	(266.157,29.403)
	(268.786,30.117)
	(271.377,30.870)
	(273.922,31.662)
	(278.839,33.371)
	(283.472,35.252)
	(287.754,37.312)
	(291.619,39.559)
	(295.000,42.000)

\path(295,42)	(298.408,46.030)
	(299.941,48.652)
	(301.329,51.485)
	(302.551,54.382)
	(303.585,57.197)
	(305.000,62.000)

\path(305,62)	(305.784,65.803)
	(306.463,70.084)
	(307.037,74.777)
	(307.286,77.255)
	(307.508,79.811)
	(307.704,82.435)
	(307.873,85.119)
	(308.017,87.854)
	(308.135,90.632)
	(308.226,93.444)
	(308.291,96.282)
	(308.330,99.137)
	(308.344,102.000)
	(308.330,104.863)
	(308.291,107.718)
	(308.226,110.556)
	(308.135,113.368)
	(308.017,116.146)
	(307.873,118.881)
	(307.704,121.565)
	(307.508,124.189)
	(307.286,126.745)
	(307.037,129.223)
	(306.463,133.916)
	(305.784,138.197)
	(305.000,142.000)

\path(305,142)	(303.585,146.803)
	(302.551,149.618)
	(301.329,152.515)
	(299.941,155.348)
	(298.408,157.970)
	(295.000,162.000)

\path(295,162)	(291.614,164.432)
	(287.737,166.654)
	(283.438,168.677)
	(278.783,170.512)
	(273.841,172.171)
	(271.284,172.938)
	(268.680,173.665)
	(266.039,174.354)
	(263.369,175.005)
	(260.678,175.621)
	(257.974,176.203)
	(255.267,176.751)
	(252.565,177.269)
	(249.876,177.756)
	(247.209,178.215)
	(244.572,178.646)
	(241.974,179.052)
	(239.424,179.433)
	(236.929,179.791)
	(232.141,180.445)
	(227.678,181.023)
	(223.608,181.538)
	(220.000,182.000)

\path(220,182)	(216.996,182.358)
	(213.796,182.669)
	(210.368,182.932)
	(206.678,183.147)
	(202.693,183.314)
	(198.381,183.434)
	(193.709,183.505)
	(191.227,183.523)
	(188.643,183.529)
	(185.953,183.523)
	(183.151,183.505)
	(180.235,183.475)
	(177.200,183.434)
	(174.041,183.380)
	(170.756,183.314)
	(167.339,183.237)
	(163.787,183.147)
	(160.095,183.045)
	(156.260,182.932)
	(152.277,182.806)
	(148.142,182.669)
	(143.851,182.520)
	(139.399,182.358)
	(134.784,182.185)
	(130.000,182.000)

\path(137.915,184.312)(130.000,182.000)(138.072,180.315)
\path(110,182)	(106.960,181.156)
	(104.030,180.327)
	(101.208,179.513)
	(98.491,178.712)
	(95.878,177.923)
	(93.364,177.145)
	(88.629,175.617)
	(84.264,174.119)
	(80.251,172.642)
	(76.569,171.176)
	(73.198,169.713)
	(70.117,168.243)
	(67.307,166.758)
	(62.418,163.703)
	(58.369,160.475)
	(55.000,157.000)

\path(55,157)	(52.271,153.338)
	(49.874,149.080)
	(48.786,146.677)
	(47.765,144.065)
	(46.804,141.224)
	(45.900,138.132)
	(45.045,134.772)
	(44.234,131.121)
	(43.463,127.160)
	(42.725,122.869)
	(42.016,118.228)
	(41.670,115.769)
	(41.329,113.216)
	(40.992,110.565)
	(40.659,107.813)
	(40.328,104.959)
	(40.000,102.000)

\path(38.876,110.169)(40.000,102.000)(42.852,109.737)
\path(40,102)	(40.212,98.783)
	(40.434,95.680)
	(40.667,92.688)
	(40.910,89.805)
	(41.165,87.028)
	(41.433,84.354)
	(41.714,81.780)
	(42.008,79.303)
	(42.641,74.632)
	(43.338,70.319)
	(44.103,66.340)
	(44.943,62.675)
	(45.863,59.302)
	(46.868,56.198)
	(47.964,53.341)
	(49.157,50.710)
	(51.853,46.036)
	(55.000,42.000)

\path(55,42)	(59.039,38.005)
	(63.708,34.259)
	(66.254,32.479)
	(68.927,30.760)
	(71.717,29.103)
	(74.614,27.506)
	(77.608,25.971)
	(80.688,24.496)
	(83.846,23.082)
	(87.070,21.729)
	(90.351,20.435)
	(93.678,19.201)
	(97.041,18.027)
	(100.431,16.913)
	(103.837,15.857)
	(107.249,14.861)
	(110.656,13.924)
	(114.050,13.045)
	(117.420,12.225)
	(120.755,11.463)
	(124.046,10.759)
	(127.282,10.113)
	(130.454,9.524)
	(133.551,8.993)
	(136.564,8.519)
	(139.481,8.102)
	(142.293,7.742)
	(144.991,7.438)
	(147.563,7.191)
	(150.000,7.000)

\path(150,7)	(153.947,6.752)
	(158.124,6.540)
	(162.514,6.369)
	(167.101,6.243)
	(171.868,6.169)
	(176.800,6.150)
	(179.323,6.163)
	(181.881,6.193)
	(184.472,6.238)
	(187.094,6.301)
	(189.745,6.382)
	(192.423,6.481)
	(195.127,6.600)
	(197.853,6.738)
	(200.600,6.896)
	(203.366,7.075)
	(206.150,7.276)
	(208.948,7.500)
	(211.759,7.746)
	(214.581,8.016)
	(217.412,8.310)
	(220.250,8.628)
	(223.093,8.972)
	(225.938,9.343)
	(228.785,9.740)
	(231.630,10.164)
	(234.472,10.617)
	(237.309,11.098)
	(240.139,11.608)
	(242.959,12.148)
	(245.769,12.719)
	(248.565,13.321)
	(251.345,13.954)
	(254.109,14.621)
	(256.853,15.320)
	(259.576,16.053)
	(262.275,16.820)
	(264.949,17.623)
	(267.596,18.461)
	(270.213,19.335)
	(272.799,20.246)
	(275.351,21.195)
	(280.348,23.207)
	(285.187,25.378)
	(289.851,27.711)
	(294.326,30.213)
	(298.594,32.887)
	(302.640,35.740)
	(306.447,38.776)
	(310.000,42.000)

\path(310,42)	(312.636,46.338)
	(313.446,49.036)
	(314.003,51.892)
	(314.372,54.762)
	(314.618,57.500)
	(315.000,62.000)

\path(315,62)	(315.326,65.989)
	(315.478,70.853)
	(315.502,73.522)
	(315.498,76.300)
	(315.471,79.150)
	(315.425,82.036)
	(315.366,84.921)
	(315.299,87.767)
	(315.228,90.540)
	(315.159,93.201)
	(315.047,98.043)
	(315.000,102.000)

\path(315,102)	(315.047,105.957)
	(315.159,110.799)
	(315.228,113.460)
	(315.299,116.233)
	(315.366,119.079)
	(315.425,121.964)
	(315.471,124.850)
	(315.498,127.700)
	(315.502,130.478)
	(315.478,133.147)
	(315.326,138.011)
	(315.000,142.000)

\path(315,142)	(314.546,146.498)
	(314.237,149.234)
	(313.811,152.102)
	(313.221,154.958)
	(312.421,157.655)
	(310.000,162.000)

\path(310,162)	(305.705,166.157)
	(300.742,170.013)
	(298.038,171.831)
	(295.200,173.576)
	(292.238,175.251)
	(289.164,176.855)
	(285.989,178.390)
	(282.723,179.857)
	(279.377,181.257)
	(275.963,182.591)
	(272.490,183.860)
	(268.971,185.065)
	(265.415,186.207)
	(261.835,187.287)
	(258.239,188.307)
	(254.641,189.266)
	(251.050,190.167)
	(247.477,191.010)
	(243.934,191.797)
	(240.430,192.527)
	(236.978,193.203)
	(233.587,193.826)
	(230.270,194.396)
	(227.036,194.914)
	(223.897,195.382)
	(220.863,195.801)
	(217.945,196.170)
	(215.155,196.493)
	(212.503,196.769)
	(210.000,197.000)

\path(210,197)	(205.542,197.229)
	(200.571,197.212)
	(197.918,197.117)
	(195.165,196.968)
	(192.322,196.767)
	(189.399,196.517)
	(186.406,196.219)
	(183.351,195.877)
	(180.244,195.493)
	(177.096,195.069)
	(173.914,194.608)
	(170.710,194.112)
	(167.493,193.584)
	(164.271,193.025)
	(161.056,192.439)
	(157.855,191.828)
	(154.680,191.195)
	(151.538,190.541)
	(148.441,189.869)
	(145.397,189.182)
	(142.416,188.482)
	(139.508,187.771)
	(136.681,187.053)
	(133.947,186.328)
	(131.314,185.601)
	(128.791,184.872)
	(124.117,183.423)
	(120.000,182.000)

\path(120,182)	(117.406,180.992)
	(114.527,179.776)
	(111.408,178.371)
	(108.097,176.797)
	(104.638,175.074)
	(101.078,173.223)
	(97.463,171.262)
	(93.838,169.213)
	(90.249,167.094)
	(86.742,164.926)
	(83.364,162.729)
	(80.159,160.522)
	(77.175,158.326)
	(74.456,156.161)
	(70.000,152.000)

\path(70,152)	(66.224,147.547)
	(64.216,144.911)
	(62.167,142.048)
	(60.106,138.989)
	(58.062,135.769)
	(56.065,132.420)
	(54.143,128.977)
	(52.326,125.471)
	(50.644,121.937)
	(49.124,118.408)
	(47.798,114.917)
	(46.693,111.498)
	(45.838,108.183)
	(45.265,105.006)
	(45.000,102.000)

\path(45,102)	(45.003,98.483)
	(45.241,94.727)
	(45.709,90.777)
	(46.400,86.677)
	(47.309,82.472)
	(48.430,78.206)
	(49.757,73.924)
	(51.284,69.670)
	(53.005,65.489)
	(54.914,61.425)
	(57.006,57.523)
	(59.275,53.827)
	(61.715,50.383)
	(64.319,47.234)
	(67.083,44.424)
	(70.000,42.000)

\path(70,42)	(72.455,40.403)
	(75.170,39.031)
	(78.108,37.873)
	(81.234,36.918)
	(84.509,36.153)
	(87.898,35.568)
	(91.365,35.149)
	(94.872,34.887)
	(98.383,34.769)
	(101.861,34.783)
	(105.270,34.919)
	(108.574,35.163)
	(111.735,35.505)
	(114.718,35.933)
	(117.485,36.435)
	(120.000,37.000)

\path(120,37)	(123.686,38.157)
	(127.789,39.871)
	(132.109,42.056)
	(136.448,44.624)
	(140.606,47.489)
	(144.383,50.565)
	(147.581,53.764)
	(150.000,57.000)

\path(150,57)	(152.100,61.606)
	(152.898,64.268)
	(153.547,67.119)
	(154.061,70.118)
	(154.456,73.229)
	(154.745,76.413)
	(154.944,79.631)
	(155.066,82.845)
	(155.126,86.017)
	(155.139,89.107)
	(155.119,92.078)
	(155.081,94.892)
	(155.038,97.509)
	(155.000,102.000)

\path(155,102)	(155.101,106.983)
	(155.185,109.886)
	(155.267,113.007)
	(155.332,116.303)
	(155.361,119.732)
	(155.339,123.251)
	(155.248,126.819)
	(155.071,130.392)
	(154.791,133.928)
	(154.391,137.386)
	(153.854,140.721)
	(153.164,143.893)
	(152.303,146.858)
	(151.254,149.575)
	(150.000,152.000)

\path(150,152)	(147.379,155.402)
	(143.715,158.655)
	(139.444,161.080)
	(135.000,162.000)

\path(135,162)	(131.816,161.889)
	(128.593,161.641)
	(125.339,161.258)
	(122.061,160.743)
	(118.767,160.099)
	(115.463,159.327)
	(112.158,158.430)
	(108.857,157.409)
	(105.569,156.268)
	(102.300,155.009)
	(99.058,153.633)
	(95.850,152.144)
	(92.684,150.543)
	(89.566,148.834)
	(86.504,147.017)
	(83.504,145.095)
	(80.576,143.072)
	(77.724,140.948)
	(74.958,138.726)
	(72.283,136.410)
	(69.708,134.000)
	(67.239,131.499)
	(64.884,128.909)
	(62.650,126.234)
	(60.544,123.474)
	(58.573,120.633)
	(56.746,117.713)
	(55.068,114.715)
	(53.547,111.643)
	(52.191,108.498)
	(51.006,105.283)
	(50.000,102.000)

\path(50,102)	(49.372,98.838)
	(49.189,95.642)
	(49.410,92.429)
	(49.994,89.217)
	(50.897,86.022)
	(52.079,82.863)
	(53.497,79.756)
	(55.110,76.718)
	(56.876,73.767)
	(58.752,70.921)
	(60.698,68.196)
	(62.671,65.610)
	(66.532,60.923)
	(70.000,57.000)

\path(70,57)	(74.687,52.765)
	(77.790,50.603)
	(81.536,48.333)
	(86.027,45.896)
	(88.585,44.597)
	(91.367,43.234)
	(94.387,41.800)
	(97.657,40.288)
	(101.190,38.691)
	(105.000,37.000)

\path(96.872,38.392)(105.000,37.000)(98.484,42.053)
\path(150,102)	(149.280,106.902)
	(149.150,109.361)
	(149.017,112.159)
	(148.869,115.221)
	(148.695,118.470)
	(148.484,121.829)
	(148.224,125.222)
	(147.904,128.572)
	(147.512,131.803)
	(147.038,134.838)
	(146.468,137.600)
	(145.000,142.000)

\path(145,142)	(142.595,146.238)
	(139.137,150.870)
	(134.861,154.817)
	(130.000,157.000)

\path(130,157)	(127.011,157.334)
	(123.947,157.269)
	(120.831,156.846)
	(117.682,156.104)
	(114.522,155.084)
	(111.372,153.825)
	(108.252,152.369)
	(105.183,150.756)
	(102.187,149.025)
	(99.285,147.216)
	(96.496,145.371)
	(93.843,143.530)
	(89.025,140.017)
	(85.000,137.000)

\path(85,137)	(81.171,134.247)
	(76.622,130.988)
	(71.725,127.243)
	(66.851,123.034)
	(62.370,118.384)
	(60.393,115.900)
	(58.654,113.314)
	(57.198,110.628)
	(56.073,107.845)
	(55.325,104.968)
	(55.000,102.000)

\path(55,102)	(55.133,98.761)
	(55.716,95.584)
	(56.701,92.478)
	(58.042,89.450)
	(59.692,86.508)
	(61.606,83.660)
	(63.737,80.912)
	(66.037,78.273)
	(70.964,73.353)
	(73.497,71.086)
	(76.014,68.960)
	(80.817,65.155)
	(85.000,62.000)

\path(85,62)	(88.659,59.444)
	(93.097,56.704)
	(95.547,55.324)
	(98.116,53.969)
	(100.782,52.662)
	(103.519,51.428)
	(106.302,50.290)
	(109.106,49.271)
	(111.908,48.396)
	(114.681,47.688)
	(117.402,47.171)
	(120.045,46.868)
	(125.000,47.000)

\path(125,47)	(127.820,47.702)
	(130.733,48.952)
	(133.648,50.644)
	(136.472,52.670)
	(139.112,54.926)
	(141.475,57.303)
	(145.000,62.000)

\path(145,62)	(146.051,64.275)
	(146.919,67.019)
	(147.620,70.145)
	(148.174,73.570)
	(148.599,77.209)
	(148.914,80.975)
	(149.136,84.785)
	(149.284,88.554)
	(149.376,92.195)
	(149.430,95.625)
	(149.466,98.759)
	(149.501,101.510)
	(149.642,105.528)
	(150.000,107.000)

\path(150,107)	(150.121,104.501)
	(150.000,102.000)

\path(475,192)	(473.100,194.697)
	(471.287,197.166)
	(467.888,201.453)
	(464.731,204.920)
	(461.742,207.625)
	(458.847,209.627)
	(455.974,210.984)
	(453.050,211.756)
	(450.000,212.000)

\path(450,212)	(445.891,211.459)
	(441.962,209.967)
	(438.264,207.719)
	(434.847,204.911)
	(431.762,201.740)
	(429.058,198.400)
	(426.788,195.088)
	(425.000,192.000)

\path(425,192)	(423.464,188.587)
	(422.195,184.770)
	(421.170,180.609)
	(420.367,176.161)
	(419.764,171.486)
	(419.339,166.643)
	(419.069,161.692)
	(418.986,159.195)
	(418.933,156.692)
	(418.909,151.701)
	(418.974,146.779)
	(419.106,141.984)
	(419.283,137.377)
	(419.483,133.016)
	(419.684,128.960)
	(419.864,125.268)
	(420.000,122.000)

\path(420,122)	(420.118,118.702)
	(420.263,114.971)
	(420.443,110.868)
	(420.667,106.455)
	(420.944,101.795)
	(421.282,96.950)
	(421.691,91.982)
	(421.925,89.471)
	(422.180,86.953)
	(422.456,84.435)
	(422.756,81.925)
	(423.430,76.962)
	(424.209,72.124)
	(425.103,67.474)
	(426.120,63.075)
	(427.269,58.988)
	(428.560,55.275)
	(430.000,52.000)

\path(430,52)	(431.363,49.305)
	(433.030,46.224)
	(435.011,42.977)
	(437.318,39.787)
	(439.959,36.873)
	(442.946,34.457)
	(446.290,32.759)
	(450.000,32.000)

\path(450,32)	(453.200,32.084)
	(456.232,32.741)
	(459.169,34.029)
	(462.085,36.008)
	(465.053,38.736)
	(468.146,42.272)
	(471.437,46.674)
	(473.180,49.218)
	(475.000,52.000)

\path(625,192)	(623.771,195.930)
	(622.606,198.751)
	(620.000,202.000)

\path(620,202)	(616.014,204.630)
	(611.786,207.102)
	(607.330,209.420)
	(602.663,211.585)
	(597.800,213.602)
	(595.301,214.556)
	(592.758,215.473)
	(590.174,216.354)
	(587.551,217.200)
	(584.891,218.011)
	(582.196,218.786)
	(579.467,219.528)
	(576.708,220.235)
	(573.919,220.909)
	(571.103,221.549)
	(568.261,222.156)
	(565.396,222.731)
	(562.510,223.273)
	(559.604,223.783)
	(556.681,224.262)
	(553.742,224.709)
	(550.790,225.125)
	(547.826,225.511)
	(544.853,225.867)
	(541.872,226.192)
	(538.885,226.489)
	(535.895,226.755)
	(532.903,226.994)
	(529.911,227.203)
	(526.921,227.385)
	(523.935,227.539)
	(520.956,227.665)
	(517.984,227.764)
	(515.023,227.837)
	(512.073,227.883)
	(509.138,227.904)
	(506.219,227.898)
	(503.317,227.868)
	(500.435,227.812)
	(497.575,227.732)
	(494.739,227.627)
	(491.929,227.499)
	(489.146,227.347)
	(486.394,227.172)
	(483.672,226.974)
	(480.985,226.753)
	(478.333,226.511)
	(475.718,226.246)
	(473.143,225.960)
	(470.610,225.653)
	(468.120,225.326)
	(463.278,224.610)
	(458.634,223.815)
	(454.202,222.944)
	(450.000,222.000)

\path(450,222)	(446.293,220.858)
	(442.183,219.144)
	(437.865,216.950)
	(433.534,214.369)
	(429.387,211.491)
	(425.618,208.409)
	(422.424,205.214)
	(420.000,202.000)

\path(420,202)	(417.942,198.171)
	(416.167,193.869)
	(414.657,189.164)
	(413.995,186.680)
	(413.391,184.122)
	(412.843,181.495)
	(412.349,178.811)
	(411.906,176.075)
	(411.512,173.298)
	(411.163,170.488)
	(410.859,167.653)
	(410.595,164.801)
	(410.370,161.941)
	(410.182,159.082)
	(410.026,156.231)
	(409.903,153.398)
	(409.808,150.591)
	(409.739,147.818)
	(409.693,145.088)
	(409.670,142.409)
	(409.664,139.790)
	(409.675,137.239)
	(409.700,134.765)
	(409.782,130.079)
	(409.888,125.802)
	(410.000,122.000)

\path(410,122)	(410.144,118.160)
	(410.361,113.838)
	(410.659,109.104)
	(410.841,106.604)
	(411.045,104.028)
	(411.272,101.384)
	(411.525,98.680)
	(411.802,95.927)
	(412.106,93.131)
	(412.436,90.303)
	(412.795,87.451)
	(413.182,84.584)
	(413.599,81.710)
	(414.047,78.838)
	(414.525,75.977)
	(415.036,73.137)
	(415.580,70.324)
	(416.158,67.549)
	(416.771,64.821)
	(417.419,62.147)
	(418.104,59.537)
	(419.587,54.542)
	(421.226,49.908)
	(423.028,45.704)
	(425.000,42.000)

\path(425,42)	(426.974,39.172)
	(429.559,36.174)
	(432.615,33.138)
	(436.000,30.195)
	(439.572,27.478)
	(443.191,25.118)
	(446.714,23.248)
	(450.000,22.000)

\path(450,22)	(454.285,20.855)
	(458.789,19.808)
	(463.495,18.862)
	(468.390,18.017)
	(470.903,17.634)
	(473.457,17.276)
	(476.051,16.945)
	(478.682,16.640)
	(481.349,16.361)
	(484.049,16.110)
	(486.781,15.885)
	(489.543,15.688)
	(492.333,15.518)
	(495.149,15.375)
	(497.989,15.261)
	(500.852,15.174)
	(503.735,15.116)
	(506.636,15.085)
	(509.554,15.084)
	(512.487,15.111)
	(515.432,15.167)
	(518.388,15.252)
	(521.353,15.367)
	(524.326,15.511)
	(527.303,15.685)
	(530.284,15.888)
	(533.266,16.122)
	(536.247,16.386)
	(539.226,16.681)
	(542.201,17.006)
	(545.170,17.362)
	(548.130,17.749)
	(551.080,18.168)
	(554.019,18.618)
	(556.943,19.100)
	(559.852,19.613)
	(562.743,20.159)
	(565.615,20.736)
	(568.465,21.347)
	(571.292,21.990)
	(574.093,22.665)
	(576.868,23.374)
	(579.613,24.116)
	(582.328,24.891)
	(585.009,25.700)
	(587.656,26.543)
	(590.267,27.420)
	(592.838,28.331)
	(595.370,29.276)
	(597.858,30.256)
	(602.702,32.321)
	(607.352,34.526)
	(611.796,36.873)
	(616.017,39.364)
	(620.000,42.000)

\path(620,42)	(622.606,45.249)
	(623.771,48.070)
	(625.000,52.000)

\path(635,192)	(633.586,195.891)
	(632.359,198.698)
	(630.000,202.000)

\path(630,202)	(625.646,205.460)
	(623.218,207.143)
	(620.637,208.794)
	(617.914,210.411)
	(615.059,211.993)
	(612.084,213.539)
	(608.999,215.047)
	(605.814,216.517)
	(602.542,217.948)
	(599.191,219.337)
	(595.774,220.685)
	(592.301,221.989)
	(588.783,223.249)
	(585.230,224.463)
	(581.654,225.631)
	(578.064,226.751)
	(574.473,227.821)
	(570.890,228.842)
	(567.326,229.811)
	(563.793,230.727)
	(560.301,231.590)
	(556.860,232.397)
	(553.482,233.149)
	(550.177,233.843)
	(546.956,234.479)
	(543.830,235.055)
	(540.810,235.570)
	(537.906,236.024)
	(535.129,236.414)
	(532.490,236.740)
	(530.000,237.000)

\path(530,237)	(525.083,237.340)
	(522.407,237.441)
	(519.598,237.496)
	(516.669,237.506)
	(513.628,237.471)
	(510.488,237.391)
	(507.258,237.267)
	(503.948,237.100)
	(500.569,236.888)
	(497.133,236.634)
	(493.648,236.337)
	(490.126,235.998)
	(486.576,235.616)
	(483.011,235.193)
	(479.439,234.729)
	(475.871,234.224)
	(472.318,233.678)
	(468.791,233.093)
	(465.300,232.467)
	(461.854,231.803)
	(458.466,231.099)
	(455.144,230.357)
	(451.900,229.576)
	(448.745,228.758)
	(445.687,227.902)
	(442.739,227.009)
	(439.911,226.079)
	(437.212,225.113)
	(434.654,224.111)
	(430.000,222.000)

\path(430,222)	(427.420,220.475)
	(424.604,218.423)
	(421.681,215.968)
	(418.781,213.234)
	(416.034,210.346)
	(413.568,207.429)
	(411.514,204.605)
	(410.000,202.000)

\path(410,202)	(408.378,198.235)
	(406.927,193.981)
	(405.640,189.308)
	(405.055,186.837)
	(404.507,184.286)
	(403.997,181.665)
	(403.522,178.983)
	(403.081,176.249)
	(402.674,173.470)
	(402.300,170.656)
	(401.958,167.815)
	(401.645,164.956)
	(401.362,162.087)
	(401.108,159.218)
	(400.881,156.357)
	(400.681,153.513)
	(400.505,150.694)
	(400.354,147.909)
	(400.226,145.167)
	(400.121,142.476)
	(400.036,139.845)
	(399.972,137.283)
	(399.927,134.798)
	(399.890,130.095)
	(399.917,125.806)
	(400.000,122.000)

\path(400,122)	(400.159,118.159)
	(400.419,113.833)
	(400.779,109.093)
	(400.998,106.590)
	(401.242,104.010)
	(401.511,101.362)
	(401.807,98.655)
	(402.128,95.897)
	(402.475,93.098)
	(402.848,90.266)
	(403.247,87.410)
	(403.673,84.540)
	(404.124,81.663)
	(404.602,78.788)
	(405.107,75.925)
	(405.638,73.082)
	(406.196,70.269)
	(406.780,67.493)
	(407.391,64.765)
	(408.030,62.092)
	(408.695,59.483)
	(410.107,54.496)
	(411.628,49.872)
	(413.259,45.683)
	(415.000,42.000)

\path(415,42)	(417.690,37.544)
	(421.257,32.624)
	(423.307,30.066)
	(425.502,27.484)
	(427.818,24.908)
	(430.230,22.367)
	(432.712,19.894)
	(435.241,17.517)
	(437.791,15.268)
	(440.338,13.177)
	(445.324,9.590)
	(450.000,7.000)

\path(450,7)	(453.861,5.486)
	(458.186,4.191)
	(462.907,3.103)
	(465.396,2.632)
	(467.958,2.209)
	(470.586,1.831)
	(473.272,1.496)
	(476.006,1.204)
	(478.782,0.952)
	(481.589,0.739)
	(484.421,0.564)
	(487.269,0.424)
	(490.124,0.319)
	(492.978,0.247)
	(495.822,0.205)
	(498.649,0.193)
	(501.450,0.209)
	(504.217,0.251)
	(506.941,0.318)
	(509.614,0.408)
	(512.228,0.520)
	(514.774,0.652)
	(517.244,0.802)
	(521.922,1.151)
	(526.197,1.554)
	(530.000,2.000)

\path(530,2)	(534.485,2.664)
	(536.930,3.081)
	(539.495,3.552)
	(542.172,4.076)
	(544.950,4.651)
	(547.818,5.275)
	(550.769,5.947)
	(553.790,6.665)
	(556.873,7.427)
	(560.007,8.232)
	(563.182,9.077)
	(566.388,9.962)
	(569.616,10.885)
	(572.856,11.843)
	(576.096,12.836)
	(579.329,13.861)
	(582.542,14.917)
	(585.727,16.002)
	(588.874,17.115)
	(591.972,18.253)
	(595.012,19.416)
	(597.983,20.601)
	(600.876,21.807)
	(603.681,23.032)
	(606.387,24.274)
	(608.985,25.532)
	(611.464,26.804)
	(616.029,29.384)
	(620.000,32.000)

\path(620,32)	(622.925,34.575)
	(626.166,38.454)
	(628.014,41.029)
	(630.074,44.106)
	(632.387,47.744)
	(635.000,52.000)

\path(625,152)	(624.765,155.731)
	(624.687,158.485)
	(625.000,162.000)

\path(625,162)	(626.913,166.609)
	(628.383,169.256)
	(630.000,172.000)
	(631.617,174.744)
	(633.087,177.391)
	(635.000,182.000)

\path(635,182)	(635.313,185.515)
	(635.235,188.269)
	(635.000,192.000)

\path(635,152)	(635.000,155.266)
	(635.000,157.000)

\path(635,157)	(635.434,159.415)
	(635.000,162.000)

\path(635,162)	(631.972,164.565)
	(628.310,166.163)
	(624.228,167.001)
	(619.941,167.285)
	(615.660,167.220)
	(611.601,167.014)
	(607.976,166.872)
	(605.000,167.000)

\path(605,167)	(602.460,167.308)
	(599.562,167.630)
	(596.357,167.961)
	(592.895,168.299)
	(589.229,168.642)
	(585.409,168.988)
	(581.486,169.333)
	(577.512,169.675)
	(573.537,170.012)
	(569.613,170.341)
	(565.791,170.659)
	(562.122,170.964)
	(558.657,171.253)
	(555.448,171.524)
	(552.545,171.774)
	(550.000,172.000)

\path(550,172)	(547.455,172.221)
	(544.552,172.455)
	(541.342,172.703)
	(537.877,172.964)
	(534.208,173.238)
	(530.386,173.524)
	(526.461,173.823)
	(522.486,174.133)
	(518.512,174.455)
	(514.589,174.788)
	(510.769,175.132)
	(507.103,175.486)
	(503.641,175.850)
	(500.437,176.224)
	(497.539,176.607)
	(495.000,177.000)

\path(495,177)	(490.522,177.328)
	(487.774,177.454)
	(484.889,177.698)
	(482.018,178.168)
	(479.312,178.969)
	(475.000,182.000)

\path(475,182)	(474.603,184.582)
	(475.000,187.000)

\path(475,187)	(475.000,188.734)
	(475.000,192.000)

\path(475,152)	(475.000,155.266)
	(475.000,157.000)

\path(475,157)	(474.603,159.418)
	(475.000,162.000)

\path(475,162)	(479.312,165.031)
	(482.018,165.832)
	(484.889,166.302)
	(487.774,166.546)
	(490.522,166.672)
	(495.000,167.000)

\path(495,167)	(497.539,167.393)
	(500.437,167.776)
	(503.641,168.150)
	(507.103,168.514)
	(510.769,168.868)
	(514.589,169.212)
	(518.512,169.545)
	(522.486,169.867)
	(526.461,170.177)
	(530.386,170.476)
	(534.208,170.762)
	(537.877,171.036)
	(541.342,171.297)
	(544.552,171.545)
	(547.455,171.779)
	(550.000,172.000)

\path(550,172)	(552.545,172.221)
	(555.448,172.455)
	(558.658,172.703)
	(562.123,172.964)
	(565.792,173.238)
	(569.614,173.524)
	(573.539,173.823)
	(577.514,174.133)
	(581.488,174.455)
	(585.411,174.788)
	(589.231,175.132)
	(592.897,175.486)
	(596.359,175.850)
	(599.563,176.224)
	(602.461,176.607)
	(605.000,177.000)

\path(605,177)	(609.478,177.328)
	(612.226,177.454)
	(615.111,177.698)
	(617.982,178.168)
	(620.688,178.969)
	(625.000,182.000)

\path(625,182)	(625.397,184.582)
	(625.000,187.000)

\path(625,187)	(625.000,188.734)
	(625.000,192.000)

\path(465,192)	(460.345,196.727)
	(456.528,199.812)
	(453.197,201.492)
	(450.000,202.000)

\path(450,202)	(445.583,201.060)
	(441.301,198.659)
	(437.618,195.429)
	(435.000,192.000)

\path(435,192)	(433.292,188.552)
	(431.906,184.707)
	(430.813,180.524)
	(429.984,176.060)
	(429.390,171.374)
	(429.002,166.526)
	(428.790,161.574)
	(428.725,156.576)
	(428.778,151.592)
	(428.920,146.679)
	(429.121,141.897)
	(429.353,137.305)
	(429.586,132.960)
	(429.791,128.922)
	(429.939,125.249)
	(430.000,122.000)

\path(430,122)	(429.984,118.977)
	(429.934,115.551)
	(429.865,111.780)
	(429.795,107.719)
	(429.740,103.425)
	(429.718,98.955)
	(429.745,94.365)
	(429.839,89.712)
	(430.016,85.053)
	(430.292,80.443)
	(430.686,75.940)
	(431.214,71.600)
	(431.892,67.479)
	(432.738,63.635)
	(433.768,60.123)
	(435.000,57.000)

\path(435,57)	(437.277,52.630)
	(440.580,47.883)
	(444.844,43.945)
	(447.314,42.649)
	(450.000,42.000)

\path(450,42)	(453.429,42.243)
	(456.838,43.835)
	(460.578,47.009)
	(465.000,52.000)

\put(315,223){\makebox(0,0)[lb]{
\raisebox{0pt}[0pt][0pt]{\shortstack[l]{{\rm\scriptsize the}}}}}
\put(315,209){\makebox(0,0)[lb]{
\raisebox{0pt}[0pt][0pt]{\shortstack[l]{{\rm\scriptsize path}}}}}
\put(315,191){\makebox(0,0)[lb]{
\raisebox{0pt}[0pt][0pt]{\shortstack[l]{\mathmode{\scriptstyle S_1S_3S_3}}}}}
\put(325,102){\makebox(0,0)[lb]{
\raisebox{0pt}[0pt][0pt]{\shortstack[l]{\mathmode{\scriptstyle S_3}}}}}
\put(245,102){\makebox(0,0)[lb]{
\raisebox{0pt}[0pt][0pt]{\shortstack[l]{\mathmode{\scriptstyle S_2}}}}}
\put(165,102){\makebox(0,0)[lb]{
\raisebox{0pt}[0pt][0pt]{\shortstack[l]{\mathmode{\scriptstyle S_1}}}}}
\put(0,48){\makebox(0,0)[lb]{
\raisebox{0pt}[0pt][0pt]{\shortstack[l]{{\rm\scriptsize the}}}}}
\put(0,34){\makebox(0,0)[lb]{
\raisebox{0pt}[0pt][0pt]{\shortstack[l]{{\rm\scriptsize path}}}}}
\put(10,16){\makebox(0,0)[lb]{
\raisebox{0pt}[0pt][0pt]{\shortstack[l]{\mathmode{\scriptstyle S_1}}}}}
\end{picture}
}}
\def\PossibleLifts
{{
\begin{picture}(380,217)(0,-10)
\put(130.000,-225.833){\arc{821.666}{4.4413}{4.9834}}
\put(170.000,-382.500){\arc{1145.000}{4.4473}{4.9775}}
\put(200.000,150.500){\arc{89.000}{3.5952}{5.8296}}
\put(240.000,78.750){\arc{242.706}{3.9926}{5.4322}}
\put(170.000,-345.833){\arc{821.666}{4.4413}{4.9834}}
\put(210.000,-502.500){\arc{1145.000}{4.4473}{4.9775}}
\put(240.000,30.500){\arc{89.000}{3.5952}{5.8296}}
\put(280.000,-41.250){\arc{242.706}{3.9926}{5.4322}}
\path(0,165)(80,165)
\path(72.000,163.000)(80.000,165.000)(72.000,167.000)
\path(140,165)(220,165)
\path(212.000,163.000)(220.000,165.000)(212.000,167.000)
\path(220,165)(300,165)
\path(292.000,163.000)(300.000,165.000)(292.000,167.000)
\path(300,165)(380,165)
\path(372.000,163.000)(380.000,165.000)(372.000,167.000)
\path(20,165)(20,170)
\path(60,165)(60,170)
\path(200,165)(200,170)
\path(240,165)(240,170)
\path(280,165)(280,170)
\path(320,165)(320,170)
\path(360,165)(360,170)
\path(160,165)(160,170)
\path(0,45)(80,45)
\path(72.000,43.000)(80.000,45.000)(72.000,47.000)
\path(140,45)(220,45)
\path(212.000,43.000)(220.000,45.000)(212.000,47.000)
\path(220,45)(300,45)
\path(292.000,43.000)(300.000,45.000)(292.000,47.000)
\path(300,45)(380,45)
\path(372.000,43.000)(380.000,45.000)(372.000,47.000)
\path(20,45)(20,50)
\path(60,45)(60,50)
\path(200,45)(200,50)
\path(240,45)(240,50)
\path(280,45)(280,50)
\path(320,45)(320,50)
\path(360,45)(360,50)
\path(160,45)(160,50)
\path(300,150)	(302.915,146.504)
	(305.000,145.000)

\path(305,145)	(308.650,144.111)
	(312.921,143.984)
	(317.607,144.365)
	(322.500,145.000)
	(327.393,145.635)
	(332.079,146.016)
	(336.350,145.889)
	(340.000,145.000)

\path(340,145)	(343.161,143.257)
	(345.000,140.000)

\path(345,140)	(343.531,136.469)
	(340.000,135.000)

\path(340,135)	(336.469,136.469)
	(335.000,140.000)

\path(335,140)	(336.839,143.257)
	(340.000,145.000)

\path(340,145)	(343.650,145.889)
	(347.921,146.016)
	(352.607,145.635)
	(357.500,145.000)
	(362.393,144.365)
	(367.079,143.984)
	(371.350,144.111)
	(375.000,145.000)

\path(375,145)	(377.085,146.504)
	(380.000,150.000)

\path(220,150)	(222.915,146.504)
	(225.000,145.000)

\path(225,145)	(228.650,144.111)
	(232.921,143.984)
	(237.607,144.365)
	(242.500,145.000)
	(247.393,145.635)
	(252.079,146.016)
	(256.350,145.889)
	(260.000,145.000)

\path(260,145)	(263.161,143.257)
	(265.000,140.000)

\path(265,140)	(263.531,136.469)
	(260.000,135.000)

\path(260,135)	(256.469,136.469)
	(255.000,140.000)

\path(255,140)	(256.839,143.257)
	(260.000,145.000)

\path(260,145)	(263.650,145.889)
	(267.921,146.016)
	(272.607,145.635)
	(277.500,145.000)
	(282.393,144.365)
	(287.079,143.984)
	(291.350,144.111)
	(295.000,145.000)

\path(295,145)	(297.085,146.504)
	(300.000,150.000)

\path(140,150)	(142.916,146.504)
	(145.000,145.000)

\path(145,145)	(148.650,144.111)
	(152.921,143.984)
	(157.607,144.365)
	(162.500,145.000)
	(167.393,145.635)
	(172.079,146.016)
	(176.350,145.889)
	(180.000,145.000)

\path(180,145)	(183.161,143.257)
	(185.000,140.000)

\path(185,140)	(183.531,136.469)
	(180.000,135.000)

\path(180,135)	(176.469,136.469)
	(175.000,140.000)

\path(175,140)	(176.839,143.257)
	(180.000,145.000)

\path(180,145)	(183.650,145.889)
	(187.921,146.016)
	(192.607,145.635)
	(197.500,145.000)
	(202.393,144.365)
	(207.079,143.984)
	(211.350,144.111)
	(215.000,145.000)

\path(215,145)	(217.085,146.504)
	(220.000,150.000)

\path(0,150)	(2.916,146.504)
	(5.000,145.000)

\path(5,145)	(8.650,144.111)
	(12.921,143.984)
	(17.607,144.365)
	(22.500,145.000)
	(27.393,145.635)
	(32.079,146.016)
	(36.350,145.889)
	(40.000,145.000)

\path(40,145)	(43.161,143.257)
	(45.000,140.000)

\path(45,140)	(43.531,136.469)
	(40.000,135.000)

\path(40,135)	(36.469,136.469)
	(35.000,140.000)

\path(35,140)	(36.839,143.257)
	(40.000,145.000)

\path(40,145)	(43.650,145.889)
	(47.921,146.016)
	(52.607,145.635)
	(57.500,145.000)
	(62.393,144.365)
	(67.079,143.984)
	(71.350,144.111)
	(75.000,145.000)

\path(75,145)	(77.084,146.504)
	(80.000,150.000)

\path(300,30)	(302.915,26.504)
	(305.000,25.000)

\path(305,25)	(308.650,24.111)
	(312.921,23.984)
	(317.607,24.365)
	(322.500,25.000)
	(327.393,25.635)
	(332.079,26.016)
	(336.350,25.889)
	(340.000,25.000)

\path(340,25)	(343.161,23.257)
	(345.000,20.000)

\path(345,20)	(343.531,16.469)
	(340.000,15.000)

\path(340,15)	(336.469,16.469)
	(335.000,20.000)

\path(335,20)	(336.839,23.257)
	(340.000,25.000)

\path(340,25)	(343.650,25.889)
	(347.921,26.016)
	(352.607,25.635)
	(357.500,25.000)
	(362.393,24.365)
	(367.079,23.984)
	(371.350,24.111)
	(375.000,25.000)

\path(375,25)	(377.085,26.504)
	(380.000,30.000)

\path(220,30)	(222.915,26.504)
	(225.000,25.000)

\path(225,25)	(228.650,24.111)
	(232.921,23.984)
	(237.607,24.365)
	(242.500,25.000)
	(247.393,25.635)
	(252.079,26.016)
	(256.350,25.889)
	(260.000,25.000)

\path(260,25)	(263.161,23.257)
	(265.000,20.000)

\path(265,20)	(263.531,16.469)
	(260.000,15.000)

\path(260,15)	(256.469,16.469)
	(255.000,20.000)

\path(255,20)	(256.839,23.257)
	(260.000,25.000)

\path(260,25)	(263.650,25.889)
	(267.921,26.016)
	(272.607,25.635)
	(277.500,25.000)
	(282.393,24.365)
	(287.079,23.984)
	(291.350,24.111)
	(295.000,25.000)

\path(295,25)	(297.085,26.504)
	(300.000,30.000)

\path(140,30)	(142.916,26.504)
	(145.000,25.000)

\path(145,25)	(148.650,24.111)
	(152.921,23.984)
	(157.607,24.365)
	(162.500,25.000)
	(167.393,25.635)
	(172.079,26.016)
	(176.350,25.889)
	(180.000,25.000)

\path(180,25)	(183.161,23.257)
	(185.000,20.000)

\path(185,20)	(183.531,16.469)
	(180.000,15.000)

\path(180,15)	(176.469,16.469)
	(175.000,20.000)

\path(175,20)	(176.839,23.257)
	(180.000,25.000)

\path(180,25)	(183.650,25.889)
	(187.921,26.016)
	(192.607,25.635)
	(197.500,25.000)
	(202.393,24.365)
	(207.079,23.984)
	(211.350,24.111)
	(215.000,25.000)

\path(215,25)	(217.085,26.504)
	(220.000,30.000)

\path(0,30)	(2.916,26.504)
	(5.000,25.000)

\path(5,25)	(8.650,24.111)
	(12.921,23.984)
	(17.607,24.365)
	(22.500,25.000)
	(27.393,25.635)
	(32.079,26.016)
	(36.350,25.889)
	(40.000,25.000)

\path(40,25)	(43.161,23.257)
	(45.000,20.000)

\path(45,20)	(43.531,16.469)
	(40.000,15.000)

\path(40,15)	(36.469,16.469)
	(35.000,20.000)

\path(35,20)	(36.839,23.257)
	(40.000,25.000)

\path(40,25)	(43.650,25.889)
	(47.921,26.016)
	(52.607,25.635)
	(57.500,25.000)
	(62.393,24.365)
	(67.079,23.984)
	(71.350,24.111)
	(75.000,25.000)

\path(75,25)	(77.084,26.504)
	(80.000,30.000)

\put(20,150){\makebox(0,0)[lb]{
\raisebox{0pt}[0pt][0pt]{\shortstack[l]{\mathmode{\scriptstyle 1}}}}}
\put(55,150){\makebox(0,0)[lb]{
\raisebox{0pt}[0pt][0pt]{\shortstack[l]{\mathmode{\scriptstyle 2}}}}}
\put(35,120){\makebox(0,0)[lb]{
\raisebox{0pt}[0pt][0pt]{\shortstack[l]{\mathmode{\scriptstyle S_1}}}}}
\put(175,120){\makebox(0,0)[lb]{
\raisebox{0pt}[0pt][0pt]{\shortstack[l]{\mathmode{\scriptstyle S_1}}}}}
\put(160,150){\makebox(0,0)[lb]{
\raisebox{0pt}[0pt][0pt]{\shortstack[l]{\mathmode{\scriptstyle 1}}}}}
\put(195,150){\makebox(0,0)[lb]{
\raisebox{0pt}[0pt][0pt]{\shortstack[l]{\mathmode{\scriptstyle 2}}}}}
\put(240,150){\makebox(0,0)[lb]{
\raisebox{0pt}[0pt][0pt]{\shortstack[l]{\mathmode{\scriptstyle 1}}}}}
\put(255,120){\makebox(0,0)[lb]{
\raisebox{0pt}[0pt][0pt]{\shortstack[l]{\mathmode{\scriptstyle S_3}}}}}
\put(275,150){\makebox(0,0)[lb]{
\raisebox{0pt}[0pt][0pt]{\shortstack[l]{\mathmode{\scriptstyle 2}}}}}
\put(320,150){\makebox(0,0)[lb]{
\raisebox{0pt}[0pt][0pt]{\shortstack[l]{\mathmode{\scriptstyle 1}}}}}
\put(335,120){\makebox(0,0)[lb]{
\raisebox{0pt}[0pt][0pt]{\shortstack[l]{\mathmode{\scriptstyle S_3}}}}}
\put(355,150){\makebox(0,0)[lb]{
\raisebox{0pt}[0pt][0pt]{\shortstack[l]{\mathmode{\scriptstyle 2}}}}}
\put(20,30){\makebox(0,0)[lb]{
\raisebox{0pt}[0pt][0pt]{\shortstack[l]{\mathmode{\scriptstyle 1}}}}}
\put(55,30){\makebox(0,0)[lb]{
\raisebox{0pt}[0pt][0pt]{\shortstack[l]{\mathmode{\scriptstyle 2}}}}}
\put(35,0){\makebox(0,0)[lb]{
\raisebox{0pt}[0pt][0pt]{\shortstack[l]{\mathmode{\scriptstyle S_1}}}}}
\put(175,0){\makebox(0,0)[lb]{
\raisebox{0pt}[0pt][0pt]{\shortstack[l]{\mathmode{\scriptstyle S_1}}}}}
\put(160,30){\makebox(0,0)[lb]{
\raisebox{0pt}[0pt][0pt]{\shortstack[l]{\mathmode{\scriptstyle 1}}}}}
\put(195,30){\makebox(0,0)[lb]{
\raisebox{0pt}[0pt][0pt]{\shortstack[l]{\mathmode{\scriptstyle 2}}}}}
\put(240,30){\makebox(0,0)[lb]{
\raisebox{0pt}[0pt][0pt]{\shortstack[l]{\mathmode{\scriptstyle 1}}}}}
\put(255,0){\makebox(0,0)[lb]{
\raisebox{0pt}[0pt][0pt]{\shortstack[l]{\mathmode{\scriptstyle S_3}}}}}
\put(275,30){\makebox(0,0)[lb]{
\raisebox{0pt}[0pt][0pt]{\shortstack[l]{\mathmode{\scriptstyle 2}}}}}
\put(320,30){\makebox(0,0)[lb]{
\raisebox{0pt}[0pt][0pt]{\shortstack[l]{\mathmode{\scriptstyle 1}}}}}
\put(335,0){\makebox(0,0)[lb]{
\raisebox{0pt}[0pt][0pt]{\shortstack[l]{\mathmode{\scriptstyle S_3}}}}}
\put(355,30){\makebox(0,0)[lb]{
\raisebox{0pt}[0pt][0pt]{\shortstack[l]{\mathmode{\scriptstyle 2}}}}}
\end{picture}
}}
\def\STU
{{
\begin{picture}(320,110)(0,-10)
\path(0,15)(0,95)
\path(2.000,87.000)(0.000,95.000)(-2.000,87.000)
\path(0,55)(30,55)(80,95)
\path(30,55)(80,15)
\path(120,15)(120,95)
\path(122.000,87.000)(120.000,95.000)(118.000,87.000)
\path(240,15)(240,95)
\path(242.000,87.000)(240.000,95.000)(238.000,87.000)
\path(200,95)(120,75)
\path(120,35)(200,15)
\path(240,75)(320,15)
\path(240,35)(320,95)
\put(85,50){\makebox(0,0)[lb]{
\raisebox{0pt}[0pt][0pt]{\shortstack[l]{\mathmode{=}}}}}
\put(210,50){\makebox(0,0)[lb]{
\raisebox{0pt}[0pt][0pt]{\shortstack[l]{\mathmode{-}}}}}
\put(15,0){\makebox(0,0)[lb]{
\raisebox{0pt}[0pt][0pt]{\shortstack[l]{\mathmode{S}}}}}
\put(135,0){\makebox(0,0)[lb]{
\raisebox{0pt}[0pt][0pt]{\shortstack[l]{\mathmode{T}}}}}
\put(255,0){\makebox(0,0)[lb]{
\raisebox{0pt}[0pt][0pt]{\shortstack[l]{\mathmode{U}}}}}
\end{picture}
}}
\def\SoFar
{{
\begin{picture}(560,158)(0,-10)
\path(0,40)(120,40)
\path(112.000,38.000)(120.000,40.000)(112.000,42.000)
\path(120,40)(240,40)
\path(232.000,38.000)(240.000,40.000)(232.000,42.000)
\path(320,40)(440,40)
\path(432.000,38.000)(440.000,40.000)(432.000,42.000)
\path(440,40)(560,40)
\path(552.000,38.000)(560.000,40.000)(552.000,42.000)
\path(100,40)(100,45)
\path(80,40)(80,45)
\path(20,40)(20,45)
\path(40,40)(40,45)
\path(220,40)(220,45)
\path(200,40)(200,45)
\path(140,40)(140,45)
\path(160,40)(160,45)
\path(420,40)(420,45)
\path(400,40)(400,45)
\path(340,40)(340,45)
\path(360,40)(360,45)
\path(540,40)(540,45)
\path(520,40)(520,45)
\path(460,40)(460,45)
\path(480,40)(480,45)
\path(520,45)	(518.153,47.823)
	(516.362,50.542)
	(514.627,53.159)
	(512.944,55.676)
	(509.732,60.421)
	(506.710,64.795)
	(503.865,68.816)
	(501.181,72.502)
	(498.644,75.873)
	(496.239,78.945)
	(493.951,81.738)
	(491.767,84.270)
	(487.648,88.623)
	(483.765,92.151)
	(480.000,95.000)

\path(480,95)	(477.463,96.645)
	(474.735,98.236)
	(471.785,99.782)
	(468.586,101.293)
	(465.106,102.778)
	(461.318,104.246)
	(457.191,105.706)
	(452.697,107.167)
	(447.807,108.639)
	(445.203,109.381)
	(442.490,110.130)
	(439.662,110.885)
	(436.718,111.650)
	(433.652,112.423)
	(430.461,113.207)
	(427.141,114.003)
	(423.690,114.812)
	(420.102,115.635)
	(416.375,116.473)
	(412.505,117.327)
	(408.489,118.199)
	(404.321,119.090)
	(400.000,120.000)

\path(408.240,120.319)(400.000,120.000)(407.421,116.404)
\path(540,45)	(537.488,48.956)
	(535.052,52.766)
	(532.690,56.434)
	(530.399,59.961)
	(528.177,63.352)
	(526.022,66.609)
	(523.930,69.737)
	(521.899,72.737)
	(519.927,75.614)
	(518.012,78.370)
	(516.150,81.009)
	(514.339,83.534)
	(510.861,88.254)
	(507.558,92.556)
	(504.409,96.465)
	(501.395,100.008)
	(498.495,103.210)
	(495.688,106.096)
	(492.956,108.693)
	(490.277,111.025)
	(487.632,113.119)
	(485.000,115.000)

\path(485,115)	(482.290,116.738)
	(479.379,118.411)
	(476.234,120.026)
	(472.826,121.593)
	(469.122,123.122)
	(465.093,124.621)
	(460.705,126.099)
	(455.929,127.567)
	(453.386,128.299)
	(450.734,129.032)
	(447.969,129.767)
	(445.088,130.504)
	(442.086,131.246)
	(438.959,131.993)
	(435.705,132.746)
	(432.318,133.507)
	(428.795,134.276)
	(425.133,135.056)
	(421.326,135.846)
	(417.372,136.648)
	(413.266,137.463)
	(409.005,138.293)
	(404.584,139.138)
	(400.000,140.000)

\path(408.231,140.499)(400.000,140.000)(407.498,136.567)
\path(460,45)	(458.215,47.785)
	(456.536,50.365)
	(453.465,54.944)
	(450.728,58.809)
	(448.266,62.035)
	(443.934,66.858)
	(440.000,70.000)

\path(440,70)	(437.320,71.529)
	(434.211,72.919)
	(430.556,74.201)
	(426.239,75.403)
	(423.796,75.983)
	(421.142,76.554)
	(418.265,77.120)
	(415.148,77.685)
	(411.778,78.251)
	(408.140,78.824)
	(404.219,79.405)
	(400.000,80.000)

\path(408.198,80.889)(400.000,80.000)(407.652,76.927)
\path(480,45)	(478.389,48.829)
	(476.859,52.380)
	(475.405,55.665)
	(474.019,58.696)
	(472.693,61.487)
	(471.420,64.051)
	(469.005,68.548)
	(466.713,72.289)
	(464.488,75.376)
	(460.000,80.000)

\path(460,80)	(456.090,82.849)
	(451.505,85.482)
	(448.905,86.735)
	(446.071,87.956)
	(442.980,89.152)
	(439.611,90.330)
	(435.942,91.498)
	(431.950,92.663)
	(427.614,93.832)
	(422.912,95.013)
	(420.416,95.610)
	(417.821,96.213)
	(415.123,96.822)
	(412.320,97.439)
	(409.409,98.065)
	(406.387,98.699)
	(403.252,99.344)
	(400.000,100.000)

\path(408.237,100.393)(400.000,100.000)(407.453,96.471)
\path(0,30)	(2.916,26.504)
	(5.000,25.000)

\path(5,25)	(7.725,24.127)
	(10.736,23.603)
	(13.990,23.379)
	(17.448,23.404)
	(21.068,23.628)
	(24.811,24.002)
	(28.635,24.476)
	(32.500,25.000)
	(36.365,25.524)
	(40.189,25.998)
	(43.932,26.372)
	(47.552,26.596)
	(51.010,26.621)
	(54.264,26.397)
	(57.275,25.873)
	(60.000,25.000)

\path(60,25)	(63.161,23.257)
	(65.000,20.000)

\path(65,20)	(63.531,16.469)
	(60.000,15.000)

\path(60,15)	(56.469,16.469)
	(55.000,20.000)

\path(55,20)	(56.839,23.257)
	(60.000,25.000)

\path(60,25)	(62.725,25.873)
	(65.736,26.397)
	(68.990,26.621)
	(72.448,26.596)
	(76.068,26.372)
	(79.811,25.998)
	(83.635,25.524)
	(87.500,25.000)
	(91.365,24.476)
	(95.189,24.002)
	(98.932,23.628)
	(102.552,23.404)
	(106.010,23.379)
	(109.264,23.603)
	(112.275,24.127)
	(115.000,25.000)

\path(115,25)	(117.084,26.504)
	(120.000,30.000)

\path(120,30)	(122.916,26.504)
	(125.000,25.000)

\path(125,25)	(127.725,24.127)
	(130.736,23.603)
	(133.990,23.379)
	(137.448,23.404)
	(141.068,23.628)
	(144.811,24.002)
	(148.635,24.476)
	(152.500,25.000)
	(156.365,25.524)
	(160.189,25.998)
	(163.932,26.372)
	(167.552,26.596)
	(171.010,26.621)
	(174.264,26.397)
	(177.275,25.873)
	(180.000,25.000)

\path(180,25)	(183.161,23.257)
	(185.000,20.000)

\path(185,20)	(183.531,16.469)
	(180.000,15.000)

\path(180,15)	(176.469,16.469)
	(175.000,20.000)

\path(175,20)	(176.839,23.257)
	(180.000,25.000)

\path(180,25)	(182.725,25.873)
	(185.736,26.397)
	(188.990,26.621)
	(192.448,26.596)
	(196.068,26.372)
	(199.811,25.998)
	(203.635,25.524)
	(207.500,25.000)
	(211.365,24.476)
	(215.189,24.002)
	(218.932,23.628)
	(222.552,23.404)
	(226.010,23.379)
	(229.264,23.603)
	(232.275,24.127)
	(235.000,25.000)

\path(235,25)	(237.085,26.504)
	(240.000,30.000)

\path(320,30)	(322.915,26.504)
	(325.000,25.000)

\path(325,25)	(327.725,24.127)
	(330.736,23.603)
	(333.990,23.379)
	(337.448,23.404)
	(341.068,23.628)
	(344.811,24.002)
	(348.635,24.476)
	(352.500,25.000)
	(356.365,25.524)
	(360.189,25.998)
	(363.932,26.372)
	(367.552,26.596)
	(371.010,26.621)
	(374.264,26.397)
	(377.275,25.873)
	(380.000,25.000)

\path(380,25)	(383.161,23.257)
	(385.000,20.000)

\path(385,20)	(383.531,16.469)
	(380.000,15.000)

\path(380,15)	(376.469,16.469)
	(375.000,20.000)

\path(375,20)	(376.839,23.257)
	(380.000,25.000)

\path(380,25)	(382.725,25.873)
	(385.736,26.397)
	(388.990,26.621)
	(392.448,26.596)
	(396.068,26.372)
	(399.811,25.998)
	(403.635,25.524)
	(407.500,25.000)
	(411.365,24.476)
	(415.189,24.002)
	(418.932,23.628)
	(422.552,23.404)
	(426.010,23.379)
	(429.264,23.603)
	(432.275,24.127)
	(435.000,25.000)

\path(435,25)	(437.085,26.504)
	(440.000,30.000)

\path(440,30)	(442.915,26.504)
	(445.000,25.000)

\path(445,25)	(447.725,24.127)
	(450.736,23.603)
	(453.990,23.379)
	(457.448,23.404)
	(461.068,23.628)
	(464.811,24.002)
	(468.635,24.476)
	(472.500,25.000)
	(476.365,25.524)
	(480.189,25.998)
	(483.932,26.372)
	(487.552,26.596)
	(491.010,26.621)
	(494.264,26.397)
	(497.275,25.873)
	(500.000,25.000)

\path(500,25)	(503.161,23.257)
	(505.000,20.000)

\path(505,20)	(503.531,16.469)
	(500.000,15.000)

\path(500,15)	(496.469,16.469)
	(495.000,20.000)

\path(495,20)	(496.839,23.257)
	(500.000,25.000)

\path(500,25)	(502.725,25.873)
	(505.736,26.397)
	(508.990,26.621)
	(512.448,26.596)
	(516.068,26.372)
	(519.811,25.998)
	(523.635,25.524)
	(527.500,25.000)
	(531.365,24.476)
	(535.189,24.002)
	(538.932,23.628)
	(542.552,23.404)
	(546.010,23.379)
	(549.264,23.603)
	(552.275,24.127)
	(555.000,25.000)

\path(555,25)	(557.084,26.504)
	(560.000,30.000)

\put(50,0){\makebox(0,0)[lb]{
\raisebox{0pt}[0pt][0pt]{\shortstack[l]{\mathmode{S_{j_m}}}}}}
\put(170,0){\makebox(0,0)[lb]{
\raisebox{0pt}[0pt][0pt]{\shortstack[l]{\mathmode{S_{j_{m-1}}}}}}}
\put(370,0){\makebox(0,0)[lb]{
\raisebox{0pt}[0pt][0pt]{\shortstack[l]{\mathmode{S_{j_1}}}}}}
\put(490,0){\makebox(0,0)[lb]{
\raisebox{0pt}[0pt][0pt]{\shortstack[l]{\mathmode{S_\nu}}}}}
\put(15,30){\makebox(0,0)[lb]{
\raisebox{0pt}[0pt][0pt]{\shortstack[l]{\mathmode{\scriptstyle 1}}}}}
\put(35,30){\makebox(0,0)[lb]{
\raisebox{0pt}[0pt][0pt]{\shortstack[l]{\mathmode{\scriptstyle 2}}}}}
\put(60,30){\makebox(0,0)[lb]{
\raisebox{0pt}[0pt][0pt]{\shortstack[l]{\mathmode{\scriptstyle m-1}}}}}
\put(50,45){\makebox(0,0)[lb]{
\raisebox{0pt}[0pt][0pt]{\shortstack[l]{\mathmode{\cdots}}}}}
\put(95,30){\makebox(0,0)[lb]{
\raisebox{0pt}[0pt][0pt]{\shortstack[l]{\mathmode{\scriptstyle m}}}}}
\put(135,30){\makebox(0,0)[lb]{
\raisebox{0pt}[0pt][0pt]{\shortstack[l]{\mathmode{\scriptstyle 1}}}}}
\put(155,30){\makebox(0,0)[lb]{
\raisebox{0pt}[0pt][0pt]{\shortstack[l]{\mathmode{\scriptstyle 2}}}}}
\put(180,30){\makebox(0,0)[lb]{
\raisebox{0pt}[0pt][0pt]{\shortstack[l]{\mathmode{\scriptstyle m-1}}}}}
\put(170,45){\makebox(0,0)[lb]{
\raisebox{0pt}[0pt][0pt]{\shortstack[l]{\mathmode{\cdots}}}}}
\put(215,30){\makebox(0,0)[lb]{
\raisebox{0pt}[0pt][0pt]{\shortstack[l]{\mathmode{\scriptstyle m}}}}}
\put(335,30){\makebox(0,0)[lb]{
\raisebox{0pt}[0pt][0pt]{\shortstack[l]{\mathmode{\scriptstyle 1}}}}}
\put(355,30){\makebox(0,0)[lb]{
\raisebox{0pt}[0pt][0pt]{\shortstack[l]{\mathmode{\scriptstyle 2}}}}}
\put(380,30){\makebox(0,0)[lb]{
\raisebox{0pt}[0pt][0pt]{\shortstack[l]{\mathmode{\scriptstyle m-1}}}}}
\put(370,45){\makebox(0,0)[lb]{
\raisebox{0pt}[0pt][0pt]{\shortstack[l]{\mathmode{\cdots}}}}}
\put(415,30){\makebox(0,0)[lb]{
\raisebox{0pt}[0pt][0pt]{\shortstack[l]{\mathmode{\scriptstyle m}}}}}
\put(455,30){\makebox(0,0)[lb]{
\raisebox{0pt}[0pt][0pt]{\shortstack[l]{\mathmode{\scriptstyle 1}}}}}
\put(475,30){\makebox(0,0)[lb]{
\raisebox{0pt}[0pt][0pt]{\shortstack[l]{\mathmode{\scriptstyle 2}}}}}
\put(500,30){\makebox(0,0)[lb]{
\raisebox{0pt}[0pt][0pt]{\shortstack[l]{\mathmode{\scriptstyle m-1}}}}}
\put(490,45){\makebox(0,0)[lb]{
\raisebox{0pt}[0pt][0pt]{\shortstack[l]{\mathmode{\cdots}}}}}
\put(535,30){\makebox(0,0)[lb]{
\raisebox{0pt}[0pt][0pt]{\shortstack[l]{\mathmode{\scriptstyle m}}}}}
\put(270,35){\makebox(0,0)[lb]{
\raisebox{0pt}[0pt][0pt]{\shortstack[l]{\mathmode{\cdots}}}}}
\put(310,130){\makebox(0,0)[lb]{
\raisebox{0pt}[0pt][0pt]{\shortstack[l]{{\rm\scriptsize to some $S_{i_m}$} \\
\mathmode{}}}}}
\put(310,110){\makebox(0,0)[lb]{
\raisebox{0pt}[0pt][0pt]{\shortstack[l]{{\rm\scriptsize to some $S_{i_{m-1}}$}
\\
\mathmode{}}}}}
\put(310,90){\makebox(0,0)[lb]{
\raisebox{0pt}[0pt][0pt]{\shortstack[l]{{\rm\scriptsize to some $S_{i_2}$} \\
\mathmode{}}}}}
\put(310,70){\makebox(0,0)[lb]{
\raisebox{0pt}[0pt][0pt]{\shortstack[l]{{\rm\scriptsize to some $S_{i_1}$} \\
\mathmode{}}}}}
\end{picture}
}}
\def\ThreeSteps
{{
\begin{picture}(516,220)(0,-10)
\path(40,33)(40,113)
\path(42.000,105.000)(40.000,113.000)(38.000,105.000)
\path(80,33)(80,113)
\path(82.000,105.000)(80.000,113.000)(78.000,105.000)
\path(120,33)(120,113)
\path(122.000,105.000)(120.000,113.000)(118.000,105.000)
\path(40,53)(120,53)
\path(80,73)(120,73)
\path(40,93)(80,93)
\path(130,103)(210,103)
\path(202.000,101.000)(210.000,103.000)(202.000,105.000)
\path(220,43)(220,58)(235,58)
	(245,68)(300,68)(300,78)
	(285,78)(275,88)(260,88)
	(260,98)(245,98)(235,108)
	(220,108)(220,123)
\path(260,43)(260,78)(275,78)
	(285,88)(300,88)(300,123)
\path(300,43)(300,58)(245,58)
	(235,68)(220,68)(220,98)
	(235,98)(245,108)(260,108)(260,123)
\path(360,103)(420,103)
\path(412.000,101.000)(420.000,103.000)(412.000,105.000)
\path(40,133)	(40.000,136.266)
	(40.000,138.000)

\path(40,138)	(39.603,140.418)
	(40.000,143.000)

\path(40,143)	(43.853,146.549)
	(48.654,148.963)
	(51.313,149.819)
	(54.093,150.480)
	(56.954,150.978)
	(59.859,151.342)
	(62.768,151.601)
	(65.643,151.785)
	(68.445,151.925)
	(71.134,152.050)
	(76.023,152.375)
	(80.000,153.000)

\path(80,153)	(83.977,153.625)
	(88.866,153.950)
	(91.555,154.075)
	(94.357,154.215)
	(97.232,154.399)
	(100.141,154.658)
	(103.046,155.022)
	(105.907,155.520)
	(108.687,156.181)
	(111.346,157.037)
	(116.147,159.451)
	(120.000,163.000)

\path(120,163)	(120.397,165.582)
	(120.000,168.000)

\path(120,168)	(120.000,169.734)
	(120.000,173.000)

\path(122.000,165.000)(120.000,173.000)(118.000,165.000)
\path(80,133)	(80.687,136.615)
	(80.917,139.330)
	(80.000,143.000)

\path(80,143)	(75.978,146.730)
	(73.272,148.154)
	(70.342,149.357)
	(67.369,150.389)
	(64.535,151.306)
	(60.000,153.000)

\path(60,153)	(55.465,154.694)
	(52.631,155.611)
	(49.658,156.643)
	(46.728,157.846)
	(44.022,159.270)
	(40.000,163.000)

\path(40,163)	(39.694,165.567)
	(40.000,168.000)

\path(40,168)	(40.000,169.734)
	(40.000,173.000)

\path(42.000,165.000)(40.000,173.000)(38.000,165.000)
\path(120,133)	(120.687,136.615)
	(120.917,139.330)
	(120.000,143.000)

\path(120,143)	(115.978,146.730)
	(113.272,148.154)
	(110.342,149.357)
	(107.369,150.389)
	(104.535,151.306)
	(100.000,153.000)

\path(100,153)	(95.465,154.694)
	(92.631,155.611)
	(89.658,156.643)
	(86.728,157.846)
	(84.022,159.270)
	(80.000,163.000)

\path(80,163)	(79.694,165.567)
	(80.000,168.000)

\path(80,168)	(80.000,169.734)
	(80.000,173.000)

\path(82.000,165.000)(80.000,173.000)(78.000,165.000)
\path(30,173)	(26.504,170.085)
	(25.000,168.000)

\path(25,168)	(24.565,164.605)
	(25.000,160.500)
	(25.435,156.395)
	(25.000,153.000)

\path(25,153)	(23.258,149.839)
	(20.000,148.000)

\path(20,148)	(16.469,149.469)
	(15.000,153.000)

\path(15,153)	(16.469,156.531)
	(20.000,158.000)

\path(20,158)	(23.258,156.161)
	(25.000,153.000)

\path(25,153)	(25.435,149.605)
	(25.000,145.500)
	(24.565,141.395)
	(25.000,138.000)

\path(25,138)	(26.504,135.915)
	(30.000,133.000)

\path(30,113)	(26.504,110.085)
	(25.000,108.000)

\path(25,108)	(24.111,104.350)
	(23.984,100.079)
	(24.365,95.393)
	(25.000,90.500)
	(25.635,85.607)
	(26.016,80.921)
	(25.889,76.650)
	(25.000,73.000)

\path(25,73)	(23.258,69.839)
	(20.000,68.000)

\path(20,68)	(16.469,69.469)
	(15.000,73.000)

\path(15,73)	(16.469,76.531)
	(20.000,78.000)

\path(20,78)	(23.258,76.161)
	(25.000,73.000)

\path(25,73)	(25.889,69.350)
	(26.016,65.079)
	(25.635,60.393)
	(25.000,55.500)
	(24.365,50.607)
	(23.984,45.921)
	(24.111,41.650)
	(25.000,38.000)

\path(25,38)	(26.504,35.915)
	(30.000,33.000)

\path(300,163)	(300.000,166.266)
	(300.000,168.000)

\path(300,168)	(299.806,170.460)
	(300.000,173.000)

\path(300,173)	(301.405,176.049)
	(303.485,179.321)
	(306.322,181.933)
	(310.000,183.000)

\path(310,183)	(313.655,181.893)
	(316.455,179.215)
	(318.528,175.929)
	(320.000,173.000)

\path(320,173)	(321.447,169.727)
	(322.754,166.020)
	(323.930,161.941)
	(324.978,157.551)
	(325.907,152.912)
	(326.722,148.086)
	(327.430,143.133)
	(327.745,140.628)
	(328.036,138.115)
	(328.304,135.601)
	(328.548,133.093)
	(328.971,128.129)
	(329.312,123.285)
	(329.576,118.621)
	(329.771,114.199)
	(329.902,110.081)
	(329.977,106.327)
	(330.000,103.000)

\path(330,103)	(329.977,99.673)
	(329.902,95.919)
	(329.771,91.801)
	(329.576,87.379)
	(329.311,82.715)
	(328.971,77.871)
	(328.548,72.907)
	(328.303,70.399)
	(328.036,67.885)
	(327.745,65.372)
	(327.430,62.867)
	(326.722,57.914)
	(325.907,53.088)
	(324.978,48.449)
	(323.929,44.059)
	(322.754,39.980)
	(321.446,36.273)
	(320.000,33.000)

\path(320,33)	(318.528,30.071)
	(316.455,26.785)
	(313.655,24.107)
	(310.000,23.000)

\path(310,23)	(306.322,24.067)
	(303.485,26.679)
	(301.405,29.951)
	(300.000,33.000)

\path(300,33)	(299.806,35.540)
	(300.000,38.000)

\path(300,38)	(300.000,39.734)
	(300.000,43.000)

\path(260,163)	(260.000,166.266)
	(260.000,168.000)

\path(260,168)	(259.661,170.426)
	(260.000,173.000)

\path(260,173)	(261.981,175.632)
	(264.253,178.128)
	(266.788,180.480)
	(269.555,182.675)
	(272.524,184.704)
	(275.665,186.555)
	(278.950,188.219)
	(282.346,189.685)
	(285.826,190.941)
	(289.359,191.979)
	(292.915,192.786)
	(296.465,193.352)
	(299.978,193.667)
	(303.425,193.721)
	(306.775,193.502)
	(310.000,193.000)

\path(310,193)	(313.261,191.886)
	(316.394,190.021)
	(319.354,187.583)
	(322.096,184.749)
	(324.576,181.697)
	(326.750,178.603)
	(328.573,175.645)
	(330.000,173.000)

\path(330,173)	(331.447,169.727)
	(332.754,166.020)
	(333.930,161.941)
	(334.978,157.551)
	(335.907,152.912)
	(336.722,148.086)
	(337.430,143.133)
	(337.745,140.628)
	(338.036,138.115)
	(338.304,135.601)
	(338.548,133.093)
	(338.971,128.129)
	(339.312,123.285)
	(339.576,118.621)
	(339.771,114.199)
	(339.902,110.081)
	(339.977,106.327)
	(340.000,103.000)

\path(340,103)	(339.977,99.673)
	(339.902,95.919)
	(339.771,91.801)
	(339.576,87.379)
	(339.311,82.715)
	(338.971,77.871)
	(338.548,72.907)
	(338.303,70.399)
	(338.036,67.885)
	(337.745,65.372)
	(337.430,62.867)
	(336.722,57.914)
	(335.907,53.088)
	(334.978,48.449)
	(333.929,44.059)
	(332.754,39.980)
	(331.446,36.273)
	(330.000,33.000)

\path(330,33)	(328.573,30.355)
	(326.751,27.397)
	(324.577,24.303)
	(322.097,21.251)
	(319.354,18.417)
	(316.394,15.979)
	(313.261,14.114)
	(310.000,13.000)

\path(310,13)	(306.775,12.498)
	(303.425,12.279)
	(299.978,12.333)
	(296.465,12.648)
	(292.915,13.214)
	(289.359,14.021)
	(285.826,15.059)
	(282.346,16.315)
	(278.950,17.781)
	(275.665,19.445)
	(272.524,21.296)
	(269.555,23.325)
	(266.788,25.520)
	(264.253,27.872)
	(261.981,30.368)
	(260.000,33.000)

\path(260,33)	(259.661,35.574)
	(260.000,38.000)

\path(260,38)	(260.000,39.734)
	(260.000,43.000)

\path(220,163)	(220.000,166.266)
	(220.000,168.000)

\path(220,168)	(219.644,170.424)
	(220.000,173.000)

\path(220,173)	(223.482,177.379)
	(227.471,181.441)
	(231.911,185.181)
	(236.745,188.597)
	(239.293,190.182)
	(241.918,191.685)
	(244.613,193.105)
	(247.371,194.441)
	(250.186,195.694)
	(253.051,196.863)
	(255.957,197.947)
	(258.899,198.945)
	(261.868,199.859)
	(264.859,200.686)
	(267.864,201.427)
	(270.876,202.081)
	(273.888,202.648)
	(276.892,203.127)
	(279.883,203.519)
	(282.852,203.821)
	(285.794,204.035)
	(288.699,204.159)
	(291.563,204.193)
	(294.377,204.137)
	(297.135,203.990)
	(299.830,203.752)
	(302.454,203.422)
	(305.000,203.000)

\path(305,203)	(307.581,202.359)
	(310.188,201.414)
	(312.805,200.194)
	(315.413,198.730)
	(317.998,197.052)
	(320.541,195.189)
	(325.436,191.029)
	(329.964,186.491)
	(333.991,181.813)
	(337.381,177.237)
	(340.000,173.000)

\path(340,173)	(341.588,169.708)
	(343.000,165.987)
	(344.247,161.898)
	(345.339,157.503)
	(346.286,152.861)
	(347.098,148.035)
	(347.785,143.085)
	(348.085,140.582)
	(348.357,138.072)
	(348.603,135.560)
	(348.824,133.056)
	(349.196,128.099)
	(349.484,123.262)
	(349.696,118.605)
	(349.844,114.189)
	(349.937,110.076)
	(349.986,106.326)
	(350.000,103.000)

\path(350,103)	(349.977,99.673)
	(349.902,95.919)
	(349.771,91.801)
	(349.576,87.379)
	(349.311,82.715)
	(348.971,77.871)
	(348.548,72.907)
	(348.303,70.399)
	(348.036,67.885)
	(347.745,65.372)
	(347.430,62.867)
	(346.722,57.914)
	(345.907,53.088)
	(344.978,48.449)
	(343.929,44.059)
	(342.754,39.980)
	(341.446,36.273)
	(340.000,33.000)

\path(340,33)	(337.841,29.045)
	(335.063,24.638)
	(331.747,20.035)
	(327.978,15.489)
	(323.835,11.258)
	(319.404,7.596)
	(314.764,4.758)
	(310.000,3.000)

\path(310,3)	(307.276,2.419)
	(304.475,1.943)
	(301.604,1.571)
	(298.671,1.303)
	(295.682,1.137)
	(292.644,1.073)
	(289.565,1.110)
	(286.451,1.248)
	(283.311,1.486)
	(280.150,1.824)
	(276.977,2.260)
	(273.798,2.794)
	(270.620,3.425)
	(267.451,4.153)
	(264.297,4.977)
	(261.167,5.896)
	(258.066,6.909)
	(255.002,8.017)
	(251.983,9.218)
	(249.015,10.511)
	(246.105,11.896)
	(243.261,13.372)
	(240.490,14.939)
	(237.799,16.595)
	(235.195,18.341)
	(232.684,20.175)
	(227.975,24.106)
	(223.728,28.383)
	(220.000,33.000)

\path(220,33)	(219.636,35.577)
	(220.000,38.000)

\path(220,38)	(220.000,39.734)
	(220.000,43.000)

\path(260,123)	(260.687,126.615)
	(260.916,129.330)
	(260.000,133.000)

\path(260,133)	(255.978,136.730)
	(253.272,138.154)
	(250.342,139.357)
	(247.369,140.389)
	(244.535,141.306)
	(240.000,143.000)

\path(240,143)	(235.465,144.694)
	(232.631,145.611)
	(229.658,146.643)
	(226.728,147.846)
	(224.022,149.270)
	(220.000,153.000)

\path(220,153)	(219.694,155.567)
	(220.000,158.000)

\path(220,158)	(220.000,159.734)
	(220.000,163.000)

\path(222.000,155.000)(220.000,163.000)(218.000,155.000)
\path(300,123)	(300.687,126.615)
	(300.916,129.330)
	(300.000,133.000)

\path(300,133)	(295.978,136.730)
	(293.272,138.154)
	(290.342,139.357)
	(287.369,140.389)
	(284.535,141.306)
	(280.000,143.000)

\path(280,143)	(275.465,144.694)
	(272.631,145.611)
	(269.658,146.643)
	(266.728,147.846)
	(264.022,149.270)
	(260.000,153.000)

\path(260,153)	(259.694,155.567)
	(260.000,158.000)

\path(260,158)	(260.000,159.734)
	(260.000,163.000)

\path(262.000,155.000)(260.000,163.000)(258.000,155.000)
\path(220,123)	(220.000,126.266)
	(220.000,128.000)

\path(220,128)	(219.603,130.418)
	(220.000,133.000)

\path(220,133)	(223.853,136.549)
	(228.654,138.963)
	(231.313,139.819)
	(234.093,140.480)
	(236.954,140.978)
	(239.859,141.342)
	(242.768,141.601)
	(245.643,141.785)
	(248.445,141.925)
	(251.134,142.050)
	(256.023,142.375)
	(260.000,143.000)

\path(260,143)	(263.977,143.625)
	(268.866,143.950)
	(271.555,144.075)
	(274.357,144.215)
	(277.232,144.399)
	(280.141,144.658)
	(283.046,145.022)
	(285.907,145.520)
	(288.687,146.181)
	(291.346,147.037)
	(296.147,149.451)
	(300.000,153.000)

\path(300,153)	(300.397,155.582)
	(300.000,158.000)

\path(300,158)	(300.000,159.734)
	(300.000,163.000)

\path(302.000,155.000)(300.000,163.000)(298.000,155.000)
\put(0,148){\makebox(0,0)[lb]{
\raisebox{0pt}[0pt][0pt]{\shortstack[l]{\mathmode{\sigma}}}}}
\put(0,68){\makebox(0,0)[lb]{
\raisebox{0pt}[0pt][0pt]{\shortstack[l]{\mathmode{D}}}}}
\put(130,113){\makebox(0,0)[lb]{
\raisebox{0pt}[0pt][0pt]{\shortstack[l]{{\rm\scriptsize steps 1 and 2}}}}}
\put(370,113){\makebox(0,0)[lb]{
\raisebox{0pt}[0pt][0pt]{\shortstack[l]{{\rm\scriptsize step 3}}}}}
\put(435,113){\makebox(0,0)[lb]{
\raisebox{0pt}[0pt][0pt]{\shortstack[l]{{\rm\scriptsize 2 components}}}}}
\put(460,93){\makebox(0,0)[lb]{
\raisebox{0pt}[0pt][0pt]{\shortstack[l]{\mathmode{N^2}}}}}
\end{picture}
}}
\def\overcross
{{
\begin{picture}(40,55)(0,-10)
\path(25,15)(40,0)
\path(15,25)(0,40)
\path(7.071,35.757)(0.000,40.000)(4.243,32.929)
\path(0,0)(40,40)
\path(35.757,32.929)(40.000,40.000)(32.929,35.757)
\end{picture}
}}
\def\selfint
{{
\begin{picture}(40,55)(0,-10)
\put(20,20){\blacken\ellipse{6}{6}}
\put(20,20){\ellipse{6}{6}}
\path(40,0)(0,40)
\path(7.071,35.757)(0.000,40.000)(4.243,32.929)
\path(0,0)(40,40)
\path(35.757,32.929)(40.000,40.000)(32.929,35.757)
\end{picture}
}}
\def\tOneTwo
{{
\begin{picture}(40,55)(0,-10)
\path(0,0)(0,40)
\path(2.000,32.000)(0.000,40.000)(-2.000,32.000)
\path(40,0)(40,40)
\path(42.000,32.000)(40.000,40.000)(38.000,32.000)
\path(0,20)(40,20)
\end{picture}
}}
\def\undercross
{{
\begin{picture}(40,55)(0,-10)
\path(40,0)(0,40)
\path(7.071,35.757)(0.000,40.000)(4.243,32.929)
\path(25,25)(40,40)
\path(35.757,32.929)(40.000,40.000)(32.929,35.757)
\path(0,0)(15,15)
\end{picture}
}}
\def\whyijk
{{
\begin{picture}(457,145)(0,-10)
\put(5,0){\makebox(0,0)[lb]{
\raisebox{0pt}[0pt][0pt]{\shortstack[l]{\mathmode{\scriptstyle i}}}}}
\put(45,0){\makebox(0,0)[lb]{
\raisebox{0pt}[0pt][0pt]{\shortstack[l]{\mathmode{\scriptstyle j}}}}}
\put(85,0){\makebox(0,0)[lb]{
\raisebox{0pt}[0pt][0pt]{\shortstack[l]{\mathmode{\scriptstyle k}}}}}
\put(125,0){\makebox(0,0)[lb]{
\raisebox{0pt}[0pt][0pt]{\shortstack[l]{\mathmode{\scriptstyle i}}}}}
\put(165,0){\makebox(0,0)[lb]{
\raisebox{0pt}[0pt][0pt]{\shortstack[l]{\mathmode{\scriptstyle j}}}}}
\put(205,0){\makebox(0,0)[lb]{
\raisebox{0pt}[0pt][0pt]{\shortstack[l]{\mathmode{\scriptstyle k}}}}}
\put(245,0){\makebox(0,0)[lb]{
\raisebox{0pt}[0pt][0pt]{\shortstack[l]{\mathmode{\scriptstyle i}}}}}
\put(285,0){\makebox(0,0)[lb]{
\raisebox{0pt}[0pt][0pt]{\shortstack[l]{\mathmode{\scriptstyle j}}}}}
\put(325,0){\makebox(0,0)[lb]{
\raisebox{0pt}[0pt][0pt]{\shortstack[l]{\mathmode{\scriptstyle k}}}}}
\put(365,0){\makebox(0,0)[lb]{
\raisebox{0pt}[0pt][0pt]{\shortstack[l]{\mathmode{\scriptstyle i}}}}}
\put(405,0){\makebox(0,0)[lb]{
\raisebox{0pt}[0pt][0pt]{\shortstack[l]{\mathmode{\scriptstyle j}}}}}
\put(445,0){\makebox(0,0)[lb]{
\raisebox{0pt}[0pt][0pt]{\shortstack[l]{\mathmode{\scriptstyle k}}}}}
\put(40,40){\blacken\ellipse{6}{6}}
\put(40,40){\ellipse{6}{6}}
\put(180,50){\blacken\ellipse{6}{6}}
\put(180,50){\ellipse{6}{6}}
\put(40,70){\blacken\ellipse{6}{6}}
\put(40,70){\ellipse{6}{6}}
\put(160,70){\blacken\ellipse{6}{6}}
\put(160,70){\ellipse{6}{6}}
\put(280,70){\blacken\ellipse{6}{6}}
\put(280,70){\ellipse{6}{6}}
\put(400,70){\blacken\ellipse{6}{6}}
\put(400,70){\ellipse{6}{6}}
\put(280,100){\blacken\ellipse{6}{6}}
\put(280,100){\ellipse{6}{6}}
\put(380,90){\blacken\ellipse{6}{6}}
\put(380,90){\ellipse{6}{6}}
\path(40,10)(40,70)(40,130)
\path(42.000,122.000)(40.000,130.000)(38.000,122.000)
\path(160,10)(160,130)
\path(162.000,122.000)(160.000,130.000)(158.000,122.000)
\path(280,10)(280,70)(280,130)
\path(282.000,122.000)(280.000,130.000)(278.000,122.000)
\path(400,10)(400,70)(400,130)
\path(402.000,122.000)(400.000,130.000)(398.000,122.000)
\path(55,55)	(51.966,58.034)
	(49.352,60.648)
	(45.204,64.797)
	(42.203,67.797)
	(40.000,70.000)

\path(40,70)	(35.550,74.358)
	(32.790,77.039)
	(29.884,79.872)
	(26.993,82.721)
	(24.276,85.450)
	(20.000,90.000)

\path(20,90)	(16.433,94.194)
	(14.281,96.810)
	(12.054,99.603)
	(9.883,102.446)
	(7.898,105.211)
	(5.000,110.000)

\path(5,110)	(3.854,112.890)
	(2.717,116.874)
	(2.115,119.423)
	(1.471,122.421)
	(0.771,125.927)
	(0.000,130.000)

\path(3.425,122.499)(0.000,130.000)(-0.508,121.769)
\path(0,10)	(-0.235,13.731)
	(-0.314,16.485)
	(0.000,20.000)

\path(0,20)	(0.712,22.400)
	(1.849,25.255)
	(5.000,30.000)

\path(5,30)	(8.390,32.136)
	(12.613,33.826)
	(17.393,35.166)
	(19.905,35.735)
	(22.454,36.253)
	(25.004,36.731)
	(27.521,37.183)
	(32.318,38.053)
	(36.570,38.960)
	(40.000,40.000)

\path(40,40)	(44.563,41.861)
	(47.339,43.068)
	(50.238,44.397)
	(53.107,45.801)
	(55.793,47.234)
	(60.000,50.000)

\path(60,50)	(63.996,53.945)
	(66.208,56.500)
	(68.415,59.272)
	(70.508,62.136)
	(72.376,64.964)
	(73.910,67.628)
	(75.000,70.000)

\path(75,70)	(76.374,74.185)
	(77.504,78.965)
	(77.982,81.633)
	(78.406,84.517)
	(78.776,87.637)
	(79.094,91.016)
	(79.362,94.675)
	(79.583,98.638)
	(79.757,102.925)
	(79.887,107.559)
	(79.936,110.012)
	(79.974,112.561)
	(80.003,115.207)
	(80.021,117.954)
	(80.030,120.804)
	(80.029,123.760)
	(80.019,126.824)
	(80.000,130.000)

\path(82.058,122.015)(80.000,130.000)(78.058,121.986)
\path(80,10)	(80.235,13.731)
	(80.314,16.485)
	(80.000,20.000)

\path(80,20)	(79.136,22.365)
	(77.748,25.160)
	(75.000,30.000)

\path(75,30)	(73.619,32.268)
	(71.648,35.289)
	(68.853,39.416)
	(67.073,42.004)
	(65.000,45.000)

\path(165,40)	(168.115,41.864)
	(170.788,43.495)
	(174.983,46.170)
	(177.936,48.260)
	(180.000,50.000)

\path(180,50)	(183.838,54.014)
	(186.043,56.572)
	(188.275,59.334)
	(190.409,62.180)
	(192.324,64.987)
	(193.895,67.634)
	(195.000,70.000)

\path(195,70)	(196.374,74.185)
	(197.504,78.965)
	(197.982,81.633)
	(198.406,84.517)
	(198.776,87.637)
	(199.094,91.016)
	(199.362,94.675)
	(199.583,98.638)
	(199.757,102.925)
	(199.887,107.559)
	(199.936,110.012)
	(199.974,112.561)
	(200.003,115.207)
	(200.021,117.954)
	(200.030,120.804)
	(200.029,123.760)
	(200.019,126.824)
	(200.000,130.000)

\path(202.058,122.015)(200.000,130.000)(198.058,121.986)
\path(200,10)	(200.235,13.731)
	(200.314,16.485)
	(200.000,20.000)

\path(200,20)	(199.146,22.368)
	(197.773,25.169)
	(195.000,30.000)

\path(195,30)	(191.900,34.575)
	(189.899,37.326)
	(187.750,40.197)
	(185.571,43.042)
	(183.476,45.720)
	(180.000,50.000)

\path(180,50)	(175.930,54.525)
	(173.355,57.260)
	(170.618,60.120)
	(167.866,62.965)
	(165.247,65.654)
	(161.000,70.000)

\path(161,70)	(156.575,74.383)
	(153.815,77.065)
	(150.906,79.894)
	(148.008,82.737)
	(145.284,85.458)
	(141.000,90.000)

\path(141,90)	(137.433,94.194)
	(135.281,96.810)
	(133.054,99.603)
	(130.883,102.446)
	(128.898,105.211)
	(126.000,110.000)

\path(126,110)	(124.854,112.890)
	(123.717,116.874)
	(123.115,119.423)
	(122.471,122.421)
	(121.771,125.927)
	(121.000,130.000)

\path(124.425,122.499)(121.000,130.000)(120.492,121.769)
\path(155,35)	(150.928,34.192)
	(147.423,33.466)
	(144.425,32.807)
	(141.876,32.200)
	(137.892,31.084)
	(135.000,30.000)

\path(135,30)	(131.366,28.366)
	(126.976,26.101)
	(122.848,23.284)
	(120.000,20.000)

\path(120,20)	(119.223,16.360)
	(119.417,13.637)
	(120.000,10.000)

\path(240,10)	(240.000,13.765)
	(240.000,16.531)
	(240.000,20.000)

\path(240,20)	(240.000,22.203)
	(240.000,25.000)
	(240.000,27.797)
	(240.000,30.000)

\path(240,30)	(239.920,34.400)
	(239.850,37.116)
	(239.787,39.982)
	(239.751,42.850)
	(239.761,45.573)
	(240.000,50.000)

\path(240,50)	(240.730,54.502)
	(241.281,57.242)
	(241.926,60.114)
	(242.643,62.971)
	(243.409,65.668)
	(245.000,70.000)

\path(245,70)	(247.771,74.852)
	(249.646,77.672)
	(251.716,80.564)
	(253.885,83.385)
	(256.053,85.994)
	(260.000,90.000)

\path(260,90)	(264.207,92.766)
	(266.893,94.199)
	(269.763,95.603)
	(272.662,96.932)
	(275.437,98.139)
	(280.000,100.000)

\path(280,100)	(283.430,101.040)
	(287.682,101.947)
	(292.479,102.817)
	(294.996,103.269)
	(297.546,103.747)
	(300.095,104.265)
	(302.607,104.834)
	(307.387,106.174)
	(311.610,107.864)
	(315.000,110.000)

\path(315,110)	(318.152,114.745)
	(319.288,117.600)
	(320.000,120.000)

\path(320,120)	(320.313,123.515)
	(320.235,126.269)
	(320.000,130.000)

\path(322.577,122.167)(320.000,130.000)(318.587,121.876)
\path(320,10)	(320.235,13.731)
	(320.313,16.485)
	(320.000,20.000)

\path(320,20)	(319.160,22.373)
	(317.809,25.181)
	(315.000,30.000)

\path(315,30)	(311.860,34.298)
	(309.925,36.739)
	(307.793,39.329)
	(305.502,42.036)
	(303.086,44.826)
	(300.582,47.666)
	(298.026,50.523)
	(295.453,53.364)
	(292.900,56.155)
	(290.402,58.864)
	(287.995,61.456)
	(283.600,66.160)
	(280.000,70.000)

\path(280,70)	(277.844,72.252)
	(274.860,75.269)
	(270.695,79.401)
	(268.061,81.995)
	(265.000,85.000)

\path(255,95)	(252.927,97.996)
	(251.147,100.585)
	(248.352,104.712)
	(246.381,107.733)
	(245.000,110.000)

\path(245,110)	(242.252,114.840)
	(240.864,117.635)
	(240.000,120.000)

\path(240,120)	(239.686,123.515)
	(239.765,126.269)
	(240.000,130.000)

\path(241.413,121.876)(240.000,130.000)(237.423,122.167)
\path(440,10)	(440.235,13.731)
	(440.313,16.485)
	(440.000,20.000)

\path(440,20)	(439.146,22.368)
	(437.773,25.169)
	(435.000,30.000)

\path(435,30)	(431.900,34.575)
	(429.899,37.326)
	(427.750,40.197)
	(425.571,43.042)
	(423.476,45.720)
	(420.000,50.000)

\path(420,50)	(415.930,54.525)
	(413.355,57.260)
	(410.618,60.120)
	(407.866,62.965)
	(405.247,65.654)
	(401.000,70.000)

\path(401,70)	(396.575,74.383)
	(393.815,77.065)
	(390.906,79.894)
	(388.008,82.737)
	(385.284,85.458)
	(381.000,90.000)

\path(381,90)	(377.433,94.194)
	(375.281,96.810)
	(373.054,99.603)
	(370.883,102.446)
	(368.898,105.211)
	(366.000,110.000)

\path(366,110)	(364.854,112.890)
	(363.717,116.874)
	(363.115,119.423)
	(362.471,122.421)
	(361.771,125.927)
	(361.000,130.000)

\path(364.425,122.499)(361.000,130.000)(360.492,121.769)
\path(360,10)	(360.000,14.045)
	(360.000,17.531)
	(360.000,20.517)
	(360.000,23.062)
	(360.000,27.062)
	(360.000,30.000)

\path(360,30)	(359.920,34.400)
	(359.850,37.116)
	(359.787,39.982)
	(359.751,42.850)
	(359.761,45.573)
	(360.000,50.000)

\path(360,50)	(360.730,54.502)
	(361.281,57.242)
	(361.926,60.114)
	(362.643,62.971)
	(363.409,65.668)
	(365.000,70.000)

\path(365,70)	(367.824,74.829)
	(369.745,77.628)
	(371.857,80.502)
	(374.049,83.313)
	(376.211,85.924)
	(380.000,90.000)

\path(380,90)	(382.064,91.740)
	(385.018,93.830)
	(389.212,96.505)
	(391.885,98.136)
	(395.000,100.000)

\path(405,105)	(409.072,105.808)
	(412.577,106.534)
	(415.575,107.193)
	(418.124,107.800)
	(422.108,108.916)
	(425.000,110.000)

\path(425,110)	(428.634,111.634)
	(433.024,113.899)
	(437.152,116.716)
	(440.000,120.000)

\path(440,120)	(440.777,123.640)
	(440.583,126.363)
	(440.000,130.000)

\path(443.431,122.501)(440.000,130.000)(439.498,121.769)
\put(95,65){\makebox(0,0)[lb]{
\raisebox{0pt}[0pt][0pt]{\shortstack[l]{\mathmode{+}}}}}
\put(215,65){\makebox(0,0)[lb]{
\raisebox{0pt}[0pt][0pt]{\shortstack[l]{\mathmode{=}}}}}
\put(335,65){\makebox(0,0)[lb]{
\raisebox{0pt}[0pt][0pt]{\shortstack[l]{\mathmode{+}}}}}
\put(30,45){\makebox(0,0)[lb]{
\raisebox{0pt}[0pt][0pt]{\shortstack[l]{\mathmode{\scriptstyle \ast}}}}}
\put(190,45){\makebox(0,0)[lb]{
\raisebox{0pt}[0pt][0pt]{\shortstack[l]{\mathmode{\scriptstyle \ast}}}}}
\put(285,90){\makebox(0,0)[lb]{
\raisebox{0pt}[0pt][0pt]{\shortstack[l]{\mathmode{\scriptstyle \ast}}}}}
\put(360,90){\makebox(0,0)[lb]{
\raisebox{0pt}[0pt][0pt]{\shortstack[l]{\mathmode{\scriptstyle \ast}}}}}
\end{picture}
}}
\def\whyijkl
{{
\begin{picture}(314,145)(0,-10)
\put(45,0){\makebox(0,0)[lb]{
\raisebox{0pt}[0pt][0pt]{\shortstack[l]{\mathmode{\scriptstyle j}}}}}
\put(5,0){\makebox(0,0)[lb]{
\raisebox{0pt}[0pt][0pt]{\shortstack[l]{\mathmode{\scriptstyle i}}}}}
\put(85,0){\makebox(0,0)[lb]{
\raisebox{0pt}[0pt][0pt]{\shortstack[l]{\mathmode{\scriptstyle k}}}}}
\put(125,0){\makebox(0,0)[lb]{
\raisebox{0pt}[0pt][0pt]{\shortstack[l]{\mathmode{\scriptstyle l}}}}}
\put(225,0){\makebox(0,0)[lb]{
\raisebox{0pt}[0pt][0pt]{\shortstack[l]{\mathmode{\scriptstyle j}}}}}
\put(185,0){\makebox(0,0)[lb]{
\raisebox{0pt}[0pt][0pt]{\shortstack[l]{\mathmode{\scriptstyle i}}}}}
\put(265,0){\makebox(0,0)[lb]{
\raisebox{0pt}[0pt][0pt]{\shortstack[l]{\mathmode{\scriptstyle k}}}}}
\put(305,0){\makebox(0,0)[lb]{
\raisebox{0pt}[0pt][0pt]{\shortstack[l]{\mathmode{\scriptstyle l}}}}}
\put(100,90){\blacken\ellipse{6}{6}}
\put(100,90){\ellipse{6}{6}}
\put(200,90){\blacken\ellipse{6}{6}}
\put(200,90){\ellipse{6}{6}}
\put(280,50){\blacken\ellipse{6}{6}}
\put(280,50){\ellipse{6}{6}}
\put(20,50){\blacken\ellipse{6}{6}}
\put(20,50){\ellipse{6}{6}}
\path(120,10)(120,20)(120,30)
	(120,50)(115,70)(100,90)
	(85,110)(80,120)(80,130)
\path(82.000,122.000)(80.000,130.000)(78.000,122.000)
\path(220,10)(220,20)(220,30)
	(220,50)(215,70)(200,90)
	(185,110)(180,120)(180,130)
\path(182.000,122.000)(180.000,130.000)(178.000,122.000)
\path(260,10)	(259.765,13.731)
	(259.687,16.485)
	(260.000,20.000)

\path(260,20)	(260.854,22.368)
	(262.227,25.169)
	(265.000,30.000)

\path(265,30)	(268.139,34.548)
	(270.174,37.277)
	(272.352,40.126)
	(274.549,42.960)
	(276.640,45.641)
	(280.000,50.000)

\path(280,50)	(283.425,54.289)
	(285.571,56.910)
	(287.819,59.688)
	(290.028,62.505)
	(292.055,65.243)
	(295.000,70.000)

\path(295,70)	(296.591,74.332)
	(297.357,77.029)
	(298.074,79.886)
	(298.719,82.758)
	(299.270,85.498)
	(300.000,90.000)

\path(300,90)	(300.239,94.427)
	(300.249,97.150)
	(300.213,100.018)
	(300.150,102.884)
	(300.080,105.600)
	(300.000,110.000)

\path(300,110)	(300.000,112.203)
	(300.000,115.000)
	(300.000,117.797)
	(300.000,120.000)

\path(300,120)	(300.000,123.469)
	(300.000,126.235)
	(300.000,130.000)

\path(302.000,122.000)(300.000,130.000)(298.000,122.000)
\path(0,10)	(-0.235,13.731)
	(-0.314,16.485)
	(0.000,20.000)

\path(0,20)	(0.854,22.368)
	(2.227,25.169)
	(5.000,30.000)

\path(5,30)	(8.139,34.548)
	(10.174,37.277)
	(12.352,40.126)
	(14.549,42.960)
	(16.640,45.641)
	(20.000,50.000)

\path(20,50)	(23.425,54.289)
	(25.571,56.910)
	(27.820,59.688)
	(30.028,62.505)
	(32.055,65.243)
	(35.000,70.000)

\path(35,70)	(36.591,74.332)
	(37.357,77.029)
	(38.074,79.886)
	(38.719,82.758)
	(39.270,85.498)
	(40.000,90.000)

\path(40,90)	(40.239,94.427)
	(40.249,97.150)
	(40.213,100.018)
	(40.150,102.884)
	(40.080,105.600)
	(40.000,110.000)

\path(40,110)	(40.000,112.203)
	(40.000,115.000)
	(40.000,117.797)
	(40.000,120.000)

\path(40,120)	(40.000,123.469)
	(40.000,126.235)
	(40.000,130.000)

\path(42.000,122.000)(40.000,130.000)(38.000,122.000)
\path(40,10)	(40.235,13.731)
	(40.314,16.485)
	(40.000,20.000)

\path(40,20)	(39.146,22.368)
	(37.773,25.169)
	(35.000,30.000)

\path(35,30)	(31.861,34.548)
	(29.826,37.277)
	(27.648,40.126)
	(25.451,42.960)
	(23.360,45.641)
	(20.000,50.000)

\path(20,50)	(16.575,54.289)
	(14.429,56.910)
	(12.180,59.688)
	(9.972,62.505)
	(7.945,65.243)
	(5.000,70.000)

\path(5,70)	(3.409,74.332)
	(2.643,77.029)
	(1.926,79.886)
	(1.281,82.758)
	(0.730,85.498)
	(0.000,90.000)

\path(0,90)	(-0.239,94.427)
	(-0.249,97.150)
	(-0.213,100.018)
	(-0.150,102.884)
	(-0.080,105.600)
	(0.000,110.000)

\path(0,110)	(0.000,112.203)
	(0.000,115.000)
	(0.000,117.797)
	(0.000,120.000)

\path(0,120)	(0.000,123.469)
	(0.000,126.235)
	(0.000,130.000)

\path(2.000,122.000)(0.000,130.000)(-2.000,122.000)
\path(80,10)	(80.000,13.034)
	(80.000,15.648)
	(80.000,19.797)
	(80.000,22.797)
	(80.000,25.000)

\path(80,25)	(80.000,27.500)
	(80.000,30.000)

\path(80,30)	(79.920,34.400)
	(79.850,37.116)
	(79.787,39.982)
	(79.751,42.850)
	(79.761,45.573)
	(80.000,50.000)

\path(80,50)	(80.730,54.502)
	(81.281,57.242)
	(81.926,60.114)
	(82.643,62.971)
	(83.409,65.668)
	(85.000,70.000)

\path(85,70)	(87.945,74.757)
	(89.972,77.495)
	(92.180,80.312)
	(94.429,83.090)
	(96.575,85.711)
	(100.000,90.000)

\path(100,90)	(103.360,94.359)
	(105.451,97.040)
	(107.648,99.874)
	(109.826,102.723)
	(111.861,105.452)
	(115.000,110.000)

\path(115,110)	(117.773,114.831)
	(119.146,117.632)
	(120.000,120.000)

\path(120,120)	(120.314,123.515)
	(120.235,126.269)
	(120.000,130.000)

\path(122.577,122.167)(120.000,130.000)(118.587,121.876)
\path(180,10)	(180.000,13.034)
	(180.000,15.648)
	(180.000,19.797)
	(180.000,22.797)
	(180.000,25.000)

\path(180,25)	(180.000,27.500)
	(180.000,30.000)

\path(180,30)	(179.920,34.400)
	(179.850,37.116)
	(179.787,39.982)
	(179.751,42.850)
	(179.761,45.573)
	(180.000,50.000)

\path(180,50)	(180.730,54.502)
	(181.281,57.242)
	(181.926,60.114)
	(182.643,62.971)
	(183.409,65.668)
	(185.000,70.000)

\path(185,70)	(187.945,74.757)
	(189.972,77.495)
	(192.180,80.312)
	(194.429,83.090)
	(196.575,85.711)
	(200.000,90.000)

\path(200,90)	(203.360,94.359)
	(205.451,97.040)
	(207.648,99.874)
	(209.826,102.723)
	(211.861,105.452)
	(215.000,110.000)

\path(215,110)	(217.773,114.831)
	(219.146,117.632)
	(220.000,120.000)

\path(220,120)	(220.314,123.515)
	(220.235,126.269)
	(220.000,130.000)

\path(222.577,122.167)(220.000,130.000)(218.587,121.876)
\path(300,10)	(300.235,13.731)
	(300.313,16.485)
	(300.000,20.000)

\path(300,20)	(299.146,22.368)
	(297.773,25.169)
	(295.000,30.000)

\path(295,30)	(291.861,34.548)
	(289.826,37.277)
	(287.648,40.126)
	(285.451,42.960)
	(283.360,45.641)
	(280.000,50.000)

\path(280,50)	(276.575,54.289)
	(274.429,56.910)
	(272.181,59.688)
	(269.972,62.505)
	(267.945,65.243)
	(265.000,70.000)

\path(265,70)	(263.409,74.332)
	(262.643,77.029)
	(261.926,79.886)
	(261.281,82.758)
	(260.730,85.498)
	(260.000,90.000)

\path(260,90)	(259.761,94.427)
	(259.751,97.150)
	(259.787,100.018)
	(259.850,102.884)
	(259.920,105.600)
	(260.000,110.000)

\path(260,110)	(260.000,112.203)
	(260.000,115.000)
	(260.000,117.797)
	(260.000,120.000)

\path(260,120)	(260.000,123.469)
	(260.000,126.235)
	(260.000,130.000)

\path(262.000,122.000)(260.000,130.000)(258.000,122.000)
\put(150,65){\makebox(0,0)[lb]{
\raisebox{0pt}[0pt][0pt]{\shortstack[l]{\mathmode{=}}}}}
\end{picture}
}}
\begin{document}
\title{Vassiliev and Quantum Invariants of Braids}

\author{Dror Bar-Natan}
\address{Department of Mathematics\\
	Harvard University\\
	Cambridge, MA 02138}
\curraddr{Institute of Mathematics\\
        The Hebrew University\\
        Giv'at-Ram, Jerusalem 91904\\
        Israel}
\email{drorbn@math.huji.ac.il}


\thanks{This work was supported by NSF grant DMS-92-03382.}

\thanks{Published in Proc.\ of Symp.\ in Appl.\ Math.\ {\bf 51} (1996)
  129--144, {\em The interface of knots and physics} (L.~H.~Kauffman,
  ed.), Amer.\ Math.\ Soc., Providence. Also available
  electronically at {\tt http://www.ma.huji.ac.il/$\sim$drorbn},
  at {\tt file://ftp.ma.huji.ac.il/drorbn}, and at \newline
  {\tt http://xxx.lanl.gov/abs/q-alg/9607001}.
}

\date{This edition: Jul.~01,~1996; \ \ First edition: November 1, 1994.}

\maketitle

\begin{abstract}
We prove that braid invariants coming from quantum $gl(N)$ separate
braids, by recalling that these invariants (properly decomposed) are
all Vassiliev invariants, showing that all Vassiliev invariants of
braids arise in this way, and reproving that Vassiliev invariants
separate braids. We discuss some corollaries of this result and of our
method of proof.
\end{abstract}

\tableofcontents
\section{Introduction} 

\subsection{The result} \lbl{TheResult}
Recall~\cite{Reshetikhin:Quasitriangle,ReshetikhinTuraev:Ribbon} that
given a list $R_1,\ldots,R_n$ of representations of the lie algebra
$gl(N)$ (or in fact, any other semi-simple Lie algebra) one can
construct an $n$-component tangle invariant (the {\em
Reshetikhin-Turaev invariant}), and in particular, an invariant
$J_{R_1,\ldots,R_n}$ of $n$-strand pure braids, with values in
$\End(R_1\otimes\cdots\otimes R_n)[[\hbar]]$, the ring of formal power
series\footnote{The Laurent polynomials in the formal parameter $q$
of~\cite{Reshetikhin:Quasitriangle,ReshetikhinTuraev:Ribbon} become
formal power series in $\hbar$ upon the substitution $q=e^\hbar$.} in
the variable $\hbar$ with coefficients in $\End(R_1\otimes\cdots\otimes
R_n)$. The main goal of this paper is to prove the following theorem

\begin{theorem} \lbl{main}
These invariants, coming from $gl(N)$ and all of its representations,
separate pure braids.
\end{theorem}

See remark~\ref{nonpurebraids} for a comment about (not necessarily pure)
braids.

The main tools we will use in proving theorem~\ref{main} are Vassiliev
invariants, chord diagrams, and weight systems. So before sketching the
proof of theorem~\ref{main} in section~\ref{sketch}, let us briefly recall
these important notions.

\subsection{Vassiliev invariants} \lbl{VassilievInvariants}
Although originally defined only for knots, the notion of `a Vassiliev
invariant' (\cite{Gusarov:New, Gusarov:nEquivalence, Vassiliev:CohKnot,
Vassiliev:Book}; see also \cite{Bar-Natan:Vassiliev, Birman:Bulletin,
BirmanLin:Vassiliev, Kontsevich:Vassiliev}) can be easily generalized
to many other classes of `knot-like' objects, such as braids, pure
braids, tangles, links, string links, knotted graphs, knots in a
3-manifold, etc. The idea is always the same. Let $\calK$ be a class
of knot-like objects --- a class of embeddings of oriented
1-dimensional objects in some 3-dimensional oriented space, perhaps
satisfying some boundary conditions, considered modulo some reasonable
notion of `isotopy'. Let $V:\calK\rightarrow\bold A$ be an invariant
(under `isotopy') with values in some Abelian group $\bold A$. It is
always possible to extend $V$ to `knot-like objects with
self-intersections' (`singular knot-like objects') by the formula
\begin{equation} \lbl{VassilievDef}
  V\left(\EEPIC{\selfint}{0.5}\right)=
    V\left(\EEPIC{\overcross}{0.5}\right)
    -V\left(\EEPIC{\undercross}{0.5}\right).
\end{equation}
(We say that the double point is `resolved into an overcrossing minus an
undercrossing').
$V$ is called {\em a Vassiliev invariant of type $m$} if its natural
extension vanishes whenever it is evaluated on an object with more than
$m$ self-intersections:
\[
  V\left(\raisebox{8pt}{$
    \underbrace{
      \EEPIC{\selfint}{0.5}\EEPIC{\selfint}{0.5}\cdots\EEPIC{\selfint}{0.5}
    }_{>m}
  $}\right) = 0.
\]

The relevance of Vassiliev invariants to our issue stems from the following
two facts:

\begin{fact} \lbl{QGareVas} (Lin~\cite{Lin:QuantumGroups}, see
also~\cite{Bar-Natan:Vassiliev,Birman:Bulletin})
The coefficient $J_{R_1,\ldots,R_n,m}$ of $\hbar^m$ in
$J_{R_1,\ldots,R_n}$ is a Vassiliev invariant of type $m$.
\end{fact}

\begin{fact} \lbl{VasSep} (Bar-Natan~\cite{Bar-Natan:Homotopy}, see
also Kohno~\cite{Kohno:deRham}) Vassiliev invariants of (pure) braids
separate (pure) braids.
\end{fact}

For the convenience of the reader, we will sketch a proof of
fact~\ref{VasSep} in section~\ref{sep}.

Differences are in many ways similar to derivatives, and as
$V\left(\EEPIC{\selfint}{0.4}\!\right)$ is defined to be a difference, one
can think of Vassiliev invariant of type $m$ as invariants whose
higher-than-$m$ derivatives vanish, or, as `polynomials of degree at
most $m$'. Just as in the case of degree $m$ polynomials, the $m$th
order `derivatives' of a type $m$ Vassiliev invariant are `constant'.
More precisely, if our knot-like object $K$ has precisely $m$
self-intersections, then by~\eqref{VassilievDef} and the definition of
Vassiliev invariants, one can replace overcrossings in $K$ by
undercrossings (and vice versa) freely, without changing the values of
any type $m$ Vassiliev invariant $V$. This means (at least in the case
where the ambient space is simply connected) that $V$ doesn't really
depend on the topology of $K$, but rather it depends only on the
combinatorial object defined by the parameter space $S$ of $K$ together
with the pairs of points in $S$ that map into each of the
self-intersections of $K$. Such pairs of points on $S$ are usually
signified by drawing a `chord' connecting them, and the resulting
combinatorial object is called {\em a chord diagram}.

The simplest and best known (see e.g.~\cite{Bar-Natan:Vassiliev,
Birman:Bulletin, BirmanLin:Vassiliev, Kontsevich:Vassiliev}) example is
that of oriented knots in oriented space. In that case, the parameter
space is an oriented circle (conventionally oriented counterclockwise
when drawn), and an example for a chord diagrams is in
figure~\ref{CDExample}.

\begin{figure}[htpb]
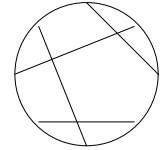

\parbox{3in}{\captionwidth=2.8in
  \caption{A chord diagram of degree $4$.} \lbl{CDExample}
}
$\qquad\EEPIC{\CDExample}{0.5}$
\end{figure}

The example that will be of interest for us is that of
$\calK=\{n$-strand pure braids$\}$. In that case, (degree $m$) chord
diagrams are diagrams made of $n$ vertical directed lines (`strands')
and $m$ horizontal lines (`chords') connecting them, as in
figure~\ref{PBExample}.

\begin{figure}[htpb]
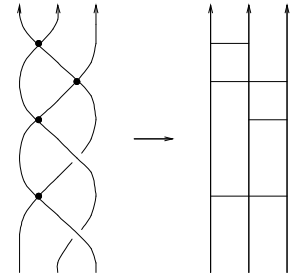

\hspace{-0.15in}
\parbox{3.2in}{\captionwidth=3.2in
  \caption{A singular 3-strand pure braid with 4 self-intersections,
    and the corresponding chord diagram.
  } \lbl{PBExample}
}$\qquad\EEPIC{\PBExample}{0.4}$
\end{figure}

Thus, to every $\bold A$-valued Vassiliev invariant $V$ of type $m$ of
pure braids, corresponds a map $W_m(V):\G{m}\calD^{pb}\rightarrow\bold
A$ defined on the $\bold Z$-module $\G{m}\calD^{pb}$ freely generated
by all degree $m$ (pure braid) chord diagrams. Notice that by stacking
chord diagrams vertically, $\calD^{pb}\eqdef\bigoplus_m\G{m}\calD^{pb}$
becomes a (non-commutative) graded algebra. As an algebra $\calD^{pb}$
is generated by the degree 1 chord diagrams $t^{ij}$, ($1\leq i\neq
j\leq n$), given by

\[ t^{ij}=t^{ji}=\EEPIC{\Omegaij}{0.5}. \]

\begin{fact} \lbl{IdealI}
$W_m(V)$ vanishes on the degree $m$ piece of the double-sided
ideal $\calI$ of $\calD^{pb}$ generated by the relations
\begin{eqnarray}
  t^{ij}t^{kl}=t^{kl}t^{ij}	& \text{when} &	|\{i,j,k,l\}|=4,
    \lbl{ijkl}	\\
  \left[t^{ik}+t^{jk},t^{ij}\right]=0	& \text{when} &	|\{i,j,k\}|=3.
    \lbl{ijk}
\end{eqnarray}
\end{fact}

Indeed, relation~\ref{ijkl} just says that double points `can be moved
across each other', as in figure~\ref{whyijkl}, while relation~\ref{ijk} is
just the $4T$ relation of~\cite{Bar-Natan:Vassiliev}, in a slightly
disguised form. Its proof is sketched in figure~\ref{whyijk}.

\begin{figure}[htpb]
\hspace{-26pt}
\parbox{3in}{\captionwidth=2.8in
  \caption{The two singular braids (which may be parts of `bigger' pure
    braids) displayed here are equivalent, implying relation~\ref{ijkl}.
  } \lbl{whyijkl}
}$\qquad\raisebox{-2mm}{$\EEPIC{\whyijkl}{0.45}$}$
\end{figure}

\begin{figure}[htpb]
\[ \EEPIC{\whyijk}{0.5} \]
\caption{The chord diagrams corresponding to the four singular braids
   displayed here are the four terms of~\eqref{ijk}. When the double points
   marked by a $\ast$ is resolved into overcrossings and undercrossings
   as in~\ref{VassilievDef}, the resulting 8 singular braids cancel in
   pairs.
} \lbl{whyijk}
\end{figure}

We set $\calA^{pb}=\calD^{pb}/\calI$, and then fact~\ref{IdealI}
just says that to every $\bold A$-valued Vassiliev invariant $V$ of
pure braids corresponds a map (a {\em weight system})
$W_m(V):\G{m}\calA^{pb}\rightarrow\bold A$. Over the rationals, the
converse is also true:

\begin{fact} \lbl{Converse}
Every weight system comes from an invariant. More precisely, if $\bold
A$ is a $\bold Q$-module and $W:\G{m}\calA^{pb}\rightarrow\bold A$ is
arbitrary, then there exists a $\bold A$-valued Vassiliev invariant $V$
of pure braids for which $W=W_m(V)$. The invariant $V$ is determined
uniquely up to invariants of lower type.
\end{fact}

Over the real numbers, fact~\ref{Converse} can be proven using the
`Kontsevich integral formula' of~\cite{Kontsevich:Vassiliev} (see
also~\cite{Bar-Natan:Vassiliev}). Over the rationals it can be deduced
from Drinfel'd's work~\cite{Drinfeld:QuasiHopf,Drinfeld:GalQQ} on
quasi-Hopf algebras. See e.g.~\cite{AltschulerFreidal:Universal,
Bar-Natan:NAT, Cartier:Construction, Kassel:QuantumGroups,
LeMurakami:Universal, Lin:Iterated, Piunikhin:Combinatorial}.

\subsection{Sketch of the proof} \lbl{sketch}

Fact~\ref{VasSep} (see section~\ref{sep})
implies that in order to prove theorem~\ref{main}, it is enough to
prove that all Vassiliev invariants of pure braids come (as in
fact~\ref{QGareVas}) from the $gl(N)$ invariants
$J_{R_1,\ldots,R_n,m}$.  Facts~\ref{IdealI} and~\ref{Converse} imply
that it is enough to do that on the level of weight systems; that is,
it is enough to prove that the weight systems corresponding to (various
traces of) the $J_{R_1,\ldots,R_n,m}$'s span the space
$\calW^{pb}=(\calA^{pb})^\ast$ of all weight systems. This is exactly
corollary~\ref{AllglN} of section~\ref{conclusion}. In
section~\ref{cors} we will prove some corollaries of theorem~\ref{main}
and discuss some related questions.

\begin{remark} The technique used in this paper appears to be a part of
a pattern --- statements about knot polynomial or quantum group
invariants of knots become simpler when restated in terms of Vassiliev
invariants, weight systems, and chord diagrams.  Chord diagrams (which
are just abstract graphs) are much more manageable objects than knots
(whose embedding into space matters), and so complicated facts about
knots become provable once stated in this simpler language. Perhaps an even
better example for this principle is the proof of the
Melvin-Morton-Rozansky conjecture in~\cite{Bar-NatanGaroufalidis:MMR}.
\end{remark}

\begin{remark} The two main problems in the theory of Vassiliev invariants
are:
\begin{itemize}
\item Do Vassiliev invariants separate knotted objects (in some class $\calK$)?
\item And do they all come from lie algebras?
\end{itemize}
In the case of $\calK=\{$knots$\}$, we know (see
e.g.~\cite{Kuperberg:Invertibility}) that the answer cannot be ``yes''
for both questions, though we don't know any of the answers.
Fact~\ref{VasSep} says that the answer to the first question is ``yes''
in the case of $\calK=\{$pure braids$\}$, and the main point in
section~\ref{paths} is to show that the answer to the second question
is also ``yes'' in that case. Another case in which both questions have an
affirmative answer is $\calK=\{$homotopy string links$\}$.
See~\cite{Bar-Natan:Homotopy}.
\end{remark}

\subsection{Acknowledgement} I would like to thank S.~Garoufalidis,
V.~Jones, D.~Kazhdan, D.~Long, P.~M.~Melvin, H.~R.~Morton, and
N.~Reshetikhin for their many useful comments, and Yael for babysitting
Assaf over the weekend in which most of this article was written.
\section{All Vassiliev invariants of pure braids come from $gl(N)$}
\lbl{paths}

In this section we will prove (as promised in section~\ref{sketch})
that the weight systems corresponding to the $gl(N)$ invariants
$J_{R_1,\ldots,R_n,m}$ span the space of all pure braid weight systems.
Denote the weight system of $J_{R_1,\ldots,R_n,m}$ by
$W_{R_1,\ldots,R_n}$ (suppressing the degree $m$, which anyway can be
read from the degree of the chord diagrams being fed into
$W_{R_1,\ldots,R_n}$). Our first step is to better understand the
$\End(R_1\otimes\cdots\otimes R_n)$-valued weight system
$W_{R_1,\ldots,R_n}$, for some fixed list $R_1,\ldots,R_n$ of
representations of $gl(N)$.

\subsection{The weight system $W_{R_1,\ldots,R_n}$} Let $R$ be the defining
($N$-dimensional) representation of $gl(N)$, and let $\{M_a\}_{a=1}^{N^2}$
be a basis of $gl(N)$, orthonormal with respect to the invariant metric
$\langle M_a,M_b\rangle=\tr_R(M_aM_b)$. Define a map
$\calD^{pb}\rightarrow\End(R_1\otimes\cdots\otimes R_n)$ by extending
the map
\begin{multline} \lbl{Casimirij}
  t^{ij}=\EEPIC{\Omegaij}{0.5}
  \\
  \mapsto\sum_{a=1}^{N^2}
  1\otimes\cdots\otimes 1\otimes\underset{i}{M_a}\otimes 1
  \cdots 1\otimes\underset{j}{M_a}\otimes 1\otimes\cdots\otimes 1
\end{multline}
to be an algebra morphism.

\begin{fact}
\begin{enumerate}
\item (Kohno~\cite{Kohno:MonRep}, see also~\cite[proposition
  2.11]{Bar-Natan:Vassiliev}) The above defined map
  $\calD^{pb}\rightarrow\End(\otimes_i R_i)$ descends to a well defined
  map $\calA^{pb}\rightarrow\End(\otimes_i R_i)$.
\item (Piunikhin~\cite{Piunikhin:Weights}, see also~\cite[remark
  4.8]{Bar-Natan:Vassiliev}) The resulting map
  $\calA^{pb}\rightarrow\End(\otimes_i R_i)$ is the weight system
  $W_{R_1,\ldots,R_n}$.
\end{enumerate}
\end{fact}

\subsection{The defining representation of $gl(N)$} Let us now
specialize to the case when $R_{1,\ldots,n}=R$, the defining
($N$-dimensional) representation of $gl(N)$.

\begin{fact} \lbl{CasimirIs} (See e.g.~\cite[exercise
6.36]{Bar-Natan:Vassiliev}) The tensor $t=\sum_a M_a\otimes
M_a\in\End(R\otimes R)$ is the `crossed identity' tensor given by the
formula
\[ t(v\otimes w)=w\otimes v. \]
\end{fact}

It is consistent with the notation in~\eqref{Casimirij} to denote the
identity operator by a directed line, and then fact~\ref{CasimirIs} becomes
the statement
\[ t=\EEPIC{\CrossedIdentity}{0.5}, \]
and~\eqref{Casimirij} becomes
\[ W_{R,\ldots,R}:\EEPIC{\Omegaij}{0.5}\mapsto\EEPIC{\Crossedij}{0.5}. \]

We can get numerical valued weight systems out of $W_{R,\ldots,R}$ by
composing it with various traces $\End(R^{\otimes n}) \rightarrow\bold
C$. I.e., if $\sigma\in\calS_n$ is a permutation of $1,\ldots,n$, it
defines in a natural way a permutation operator (also denoted by the
letter $\sigma$) in $\End(R^{\otimes n})$, and then one can define
$\tr_\sigma:\End(R^{\otimes n})\rightarrow\bold C$ by $\tr_\sigma E=\tr
\sigma\circ E$. Let $W_\sigma$ denote the numerical valued weight
system $\tr_\sigma\circ W_{R,\ldots,R}$.

\begin{proposition} \lbl{DefiningAlgo} Let $D\in\calD^{pb}$ be a chord
diagram. Then $W_\sigma(D)$ can be computed in three steps as below
(and as in figure~\ref{ThreeSteps}):
\begin{enumerate}
\item Append the permutation $\sigma$ to the top of $D$, and `close around'.
\item Replace all chords in $D$ by `crossings'.
\item Count the number $c$ of components in the resulting diagram.
  $W_\sigma(D)$ is then $N^c$.
\end{enumerate}
\end{proposition}

\begin{figure}[htpb]
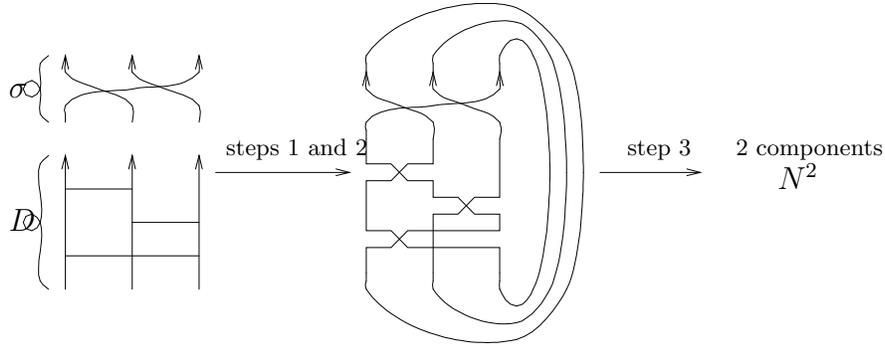

\[ \EEPIC{\ThreeSteps}{0.7} \]
\caption{The three steps of computing $W_\sigma(D)$.}
\lbl{ThreeSteps} \end{figure}

\subsection{Tensor powers of the defining representation} \lbl{TensorPowers}
The following definition and fact allow us to compute
$W_{R_1,\ldots,R_n}(D)$ whenever $D$ is a chord diagram and
$\left\{R_i=R^{\otimes k_i}\right\}$ are tensor powers of the defining
representation.

\begin{definition} Let $k=(k_1,\ldots,k_n)$ be a sequence of
non negative integers, let $|k|=\sum_ik_i$ be their sum, and let
$l_i=\sum_{\alpha<i}k_\alpha$. Define an algebra morphism
$\Delta^k:\calD^{pb}_n\rightarrow\calD^{pb}_{|k|}$ by setting
\[ \Delta^k(t^{ij})=
  \sum_{i'=l_i+1}^{l_{i+1}}
  \sum_{j'=l_j+1}^{l_{j+1}}
  t^{i'j'}.
\]
This map descends to a morphism (called by the same name)
$\Delta^k:\calA^{pb}_n\rightarrow\calA^{pb}_{|k|}$.
\end{definition}

$\Delta^k$ can be thought of as replacing the $i$th strand in $D$
by a bundle of $k_i$ strands (for each $i$) and summing over all
possible ways of `lifting' the chords from their original strand to the
bundle that replaced it. For example,
$\Delta^{(22)}(t^{12})=t^{13}+t^{14}+t^{23}+t^{24}$, or
\[ \Delta^{(22)}\left(\EEPIC{\tOneTwo}{0.5}\right)=\EEPIC{\DeltaExample}{0.5}.
\]

\begin{fact} \lbl{CommDiag} (See proposition 6.18
of~\cite{Bar-Natan:Vassiliev}) The following diagram commutes:
\[ {
  \def\a{\calA^{pb}_n}
  \def\b{\calA^{pb}_{|k|}}
  \def\c{\End(R^{\otimes|k|})}
  \def\ab{\Delta^k}
  \def\ac{W_{R^{\otimes k_1},\ldots,R^{\otimes k_n}}}
  \def\bc{W_{R,\ldots,R\ (|k|\text{ times})}}
  \EEPIC{\CommDiag}{0.75}
} \]
\end{fact}

Given a list $k=(k_1,\ldots,k_n)$ and a permutation
$\sigma\in\calS_{|k|}$ of the integers $1$ through $|k|$, we consider
the numerical valued weight system
\[ W_{k,\sigma}
  =\tr_\sigma\circ W_{R^{\otimes k_1},\ldots,R^{\otimes k_n}}.
\]
Fact~\ref{CommDiag} and
proposition~\ref{DefiningAlgo} show that we can compute
$W_{k,\sigma}(D)$ for any chord diagram $D$ following the steps below:
\begin{enumerate}
\item For each $i$, replace the $i$th strand in $D$ by a bundle of
  $k_i$ strands and consider all possible ways $\{D_\alpha\}$ of
  lifting all chords.
\item For each $\alpha$, append the permutation $\sigma$ at the top of
  $D_\alpha$, close around, replace all chords by crossings, and call
  the result $D'_\alpha$.
\item $W_{k,\sigma}(D)$ is $\sum_\alpha N^{c_\alpha}$, where $c_\alpha$ is
  the number of components in $D'_\alpha$.
\end{enumerate}

\subsection{Paths} \lbl{Paths} It is easy to see that $W_{k,\sigma}(D)$
does not change if $\sigma\in\calS_{|k|}$ is conjugated by a
permutation in $\calS_{k_1}\times\cdots\times \calS_{k_n}\subset
\calS_{|k|}$, as this just corresponds to permuting the strands within
each bundle. A {\em path}, defined below, is basically a pair
$(k,\sigma)$, with the redundancy in $\sigma$ removed.

\begin{definition} A connected path is a word in the $n$ letters
$S_1$ through $S_n$. A path is an unordered list (possibly with
multiplicities) of connected paths.
\end{definition}

Figure~\ref{PairPath} explains by an example how a pair $(k,\sigma)$ as
above determines a path $P$ (up to cyclically permuting the letters of each
word in $P$), and how a path $P$ determines a pair $(k,\sigma)$ (up to
conjugating $\sigma$ by a permutation in
$\calS_{k_1}\times\cdots\times \calS_{k_n}$).

\begin{figure}[htpb]
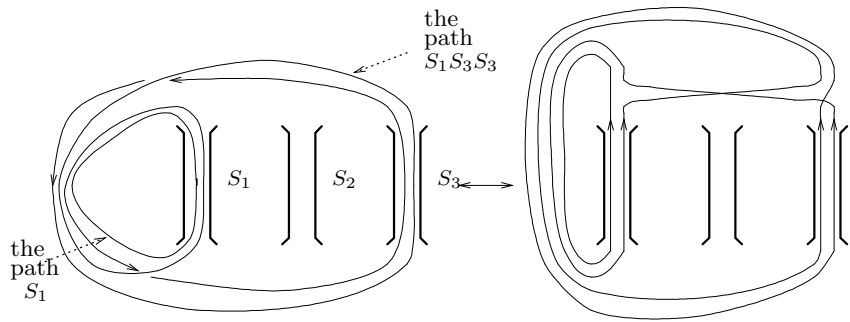

\[ \EEPIC{\PairPath}{0.55} \]
\[ \parbox{4.5in}{
  \captionwidth=4.5in
  \caption{
  The path $P=(S_1,S_1S_3S_3)$ and the pair
    $(k=(2,0,2),$ $\sigma=(1)(234))$ correspond.
  }
} \]
\lbl{PairPath} \end{figure}

We thus find that to every path $P$ corresponds a weight system
$W_P=W_{k,\sigma}$, where $(k,\sigma)$ correspond to $P$ as in
figure~\ref{PairPath}. The algorithm in section~\ref{TensorPowers} becomes
the following algorithm for computing $W_P(D)$, where $D$ is a chord
diagram of degree $m$:
\begin{enumerate}
\item For each connected component of $P$ write a long interval,
  subdivided into shorter subintervals corresponding to the letters
  making up that component. Mark sites corresponding to the integers
  $1,\ldots,m$ in order along each of the subintervals. For example, if
  $m=2$, the path $P=(S_1,S_1S_3S_3)$ becomes:
  \begin{equation} \lbl{LinearPath}
    \EEPIC{\LinearPath}{0.5}
  \end{equation}
\item Consider all liftings $\{D_\alpha\}$ of $D$ to the picture just
  drawn, where each end of each chord is lifted to one of the sites in
  the picture, so that if the (say) 7th chord in $D$ is $t^{23}$, then
  its ends are lifted to sites on marked by the integer 7 on
  subintervals corresponding to the letters $S_2$ and $S_3$. For
  example, the chord diagram $D=t^{13}t^{13}$ has 16 possible liftings
  to the path $(S_1,S_1S_3S_3)$, as the first $t^{13}$ in $D$ can be
  lifted to become either one of the 4 chords in the first figure
  below, and the second $t^{13}$ in $D$ can be lifted to become either
  one of the 4 chords in the second figure:
  \[ \EEPIC{\PossibleLifts}{0.5} \]
  As a second example, notice that the chord diagram $t^{13}t^{23}$ has
  no liftings to the path $(S_1,S_1S_3S_3)$, as that path does not pass
  through the strand $S_2$ at all, and there's nowhere to lift $t^{23}$
  to.
\item Replace all chords in each $D_\alpha$ by `bridges' (notice that
  the `crossings' of section~\ref{TensorPowers} become `bridges' in the
  current picture), erase all markings, close all intervals into loops,
  and call the result $D'_\alpha$:
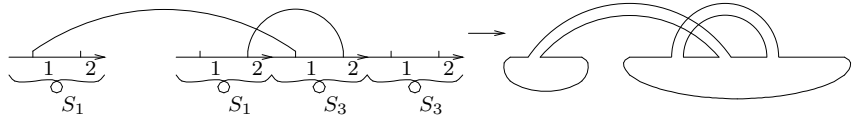
  \[ \EEPIC{\Bridged}{0.5} \]
\item $W_P(D)$ is $\sum_\alpha N^{c_\alpha}$, where $c_\alpha$ is
  the number of components in $D'_\alpha$.
\end{enumerate}

\subsection{The conclusion of the proof} \lbl{conclusion}

Define the {\em order} of a chord $t^{ij}$ to be $\max(i,j)$. We say
that a chord diagram is {\em non-decreasing} if its chords appear in
non-decreasing order. We say that it is {\em flat} if all the chords in
it are of the same order. For example, $t^{13}t^{23}t^{13}$ is flat and
non-decreasing, $t^{12}t^{23}t^{13}$ is non-decreasing but not flat,
and $t^{13}t^{12}$ is neither flat nor non-decreasing. The following
fact is a trivial consequence of the relations~\eqref{ijk}
and~\eqref{ijkl} defining $\calA^{pb}$.

\begin{fact} \lbl{NonDecGen}
(Drinfel'd~\cite{Drinfeld:GalQQ}) $\calA^{pb}$ is
generated by non-decreasing chord diagrams.
\end{fact}

The following two propositions make the key technical observation of
this paper.

\begin{proposition} \lbl{TechLemma} The weight systems $W_P$ separate
flat chord diagrams. More precisely, if
$D_\nu(i_1,\ldots,i_m)=t^{i_1\nu}\cdots t^{i_m\nu}$ is a flat chord
diagram of degree $m$ (that is, $i_\alpha<\nu$ for all $\alpha$) and
$P_\nu(j_1,\ldots,j_m)$ (with $j_\alpha<\nu$ for all $\alpha$) is the
connected path $S_{j_m}S_{j_{m-1}}\ldots S_{j_1}S_\nu$, and if
$C_{N,p}$ is the operator that maps a polynomial in $N$ to the
coefficient of $N^p$ in it, then
\[
  C_{N,m+1}W_{P_\nu(j_1,\ldots,j_m)}(D_\nu(i_1,\ldots,i_m))=
  \begin{cases}
    1 & \text{if }i_\alpha=j_\alpha\text{ for all }\alpha, \\
    0 & \text{otherwise.}
  \end{cases}
\]
\end{proposition}

\begin{pf}
Consider the diagrams $D_\alpha$ and the corresponding numbers
$c_\alpha$ that appear when
$W_{P_\nu(j_1,\ldots,j_m)}(D_\nu(i_1,\ldots,i_m))$ is computed using
the algorithm of section~\ref{Paths}. It is easy to show that the
maximal possible value for the $c_\alpha$'s is $m+1$, and that this
maximum is attained iff the corresponding diagram $D_\alpha$ has no
chord intersections:
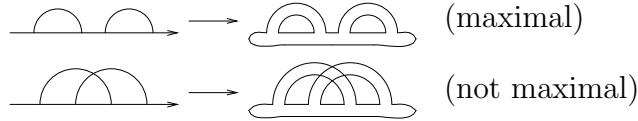
\[ \EEPIC{\MaxExample}{0.5} \]
Therefore, we only care about those $D_\alpha$'s in which there are no
chord intersections. Notice that the right ends of the chords in
$D_\nu(i_1,\ldots,i_m)$ can only be lifted (in order) to the last
interval $(S_\nu)$ in $P_\nu(j_1,\ldots,j_m)$ --- they simply have
nowhere else to go.  The left ends of these chords should be lifted to
one of the earlier intervals in $P_\nu(j_1,\ldots,j_m)$, and so far the
picture is:
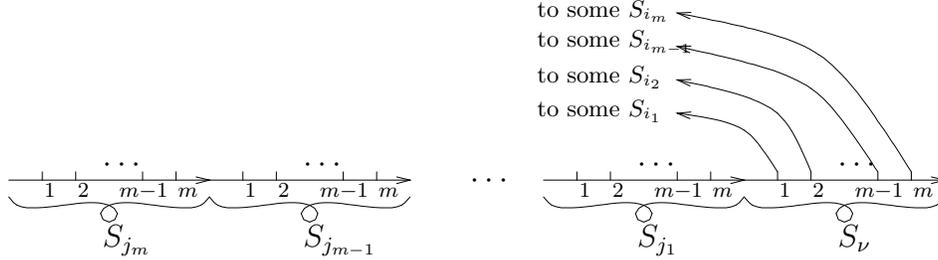
\[ \EEPIC{\SoFar}{0.7} \]
It is only possible to connect the left ends of the chords above subject to
the restrictions and without introducing chord intersections if
$i_\alpha=j_\alpha$ for all $1\leq\alpha\leq m$, and this is then possible
in a unique way.
\end{pf}

In order to show that the $W_P$'s separate non-decreasing chord
diagrams, we have to work just a little bit harder. By mapping
$D=D_\nu(i_1,\ldots,i_m)$ to $P(D)=P_\nu(i_1,\ldots,i_m)$ we see that
to every flat chord diagram $D$ corresponds a connected path $P(D)$.
This map can be extended to non-decreasing chord diagrams --- simply
write any non-decreasing chord diagram $D$ as an increasing product of
flat chord diagrams, and let $P(D)$ be the list of connected paths
corresponding to the flat parts of $D$ (including a connected path
$S_\nu$ for each $1\leq\nu\leq n$ for which the order $\nu$ part of $D$
is empty). Let the profile of a diagram $D$ be the sequence
$(l_n,l_{n-1},\ldots,l_1)$, where $l_\nu$ is the number of chords in
$D$ which are of order $\nu$. Let us order all degree $m$ chord
diagrams using some order $\prec$ that refines the lexicographic order
on their profiles. The following proposition is a generalization of
proposition~\ref{TechLemma}, and its proof is very similar.

\begin{proposition} \lbl{MainTech}
The $W_P$'s separate non-decreasing chord diagrams.  More precisely, if
$D_1$ and $D_2$ are chord diagrams of degree $m$, then
\[ C_{N,m+n}W_{P(D_1)}(D_2) = \begin{cases}
  1 & \text{if }D_1=D_2, \\
  0 & \text{if }D_1\prec D_2.
\end{cases} \]
\end{proposition}

\begin{pf}
Consider the diagrams $D_\alpha$ and the numbers $c_\alpha$ that appear
when $W_{P(D_1)}(D_2)$ is computed using the algorithm of
section~\ref{Paths}. It is easy to show that the maximal possible value
for the $c_\alpha$'s is $m+n$, and that this maximum is attained iff
the corresponding diagram $D_\alpha$ has no chord intersections, and no
chords connecting different components of $P(D_1)$. Therefore, we only
care about those $D_\alpha$'s in which there are no chord intersections
and no such connecting chords. If the profile of $D_1$ is
lexicographically smaller than the profile of $D_2$, there is simply no
room to fit a lifting of $D_2$ on $P(D_1)$ without chord intersections
or connecting chords. If the profiles of $D_1$ and $D_2$ are the same,
the flat components of $D_2$ must be lifted to the corresponding
connected paths in $P(D_1)$, and the same argument as in the proof of
proposition~\ref{TechLemma} shows that such a lifting (having no chord
intersections or connecting chords) exists iff $D_1=D_2$, and in that
case it is unique.
\end{pf}

\begin{corollary} \lbl{AllglN}
The weight systems corresponding to (traces of) the $gl(N)$
invariants $J_{R_1,\ldots,R_n,m}$ span the space
$\calW^{pb}=(\calA^{pb})^\ast$ of all weight systems. \qed
\end{corollary}
\section{Vassiliev invariants separate pure braids} \lbl{sep}

Just for the amusement of the reader, we include here the `moral reason'
for why Vassiliev invariants separate pure braids. As proofs of this
fact are available elsewhere (see e.g.~\cite{Bar-Natan:Homotopy,
Kohno:deRham}), we will leave the details of the proof presented here as an
exercise to the reader.

\begin{fact}
\begin{enumerate}
\item Every pure braid has a unique presentation as a `combed
  braid', as in figure~\ref{CombedBraid}. Equivalently, the group $K_n$ of
  pure braids on $n$ strands is a semi-direct product of free groups:
  \begin{equation} \lbl{SemiDirectF}
    K_n=F_1\ltimes F_2\ltimes\cdots\ltimes F_{n-1}.
  \end{equation}
\item (Kohno~\cite{Kohno:MonRep}) The non-decreasing chord diagrams form
  a basis of $\calA^{pb}$.  Equivalently, the algebra $\calA^{pb}_n$ of
  chord diagrams on $n$ strands is a semi-direct product of free
  associative algebras:
  \begin{equation} \lbl{SemiDirectFA}
    \calA^{pb}_n=FA_1\ltimes FA_2\ltimes\cdots\ltimes FA_{n-1}.
  \end{equation}
\end{enumerate}
\end{fact}

\begin{figure}[htpb]
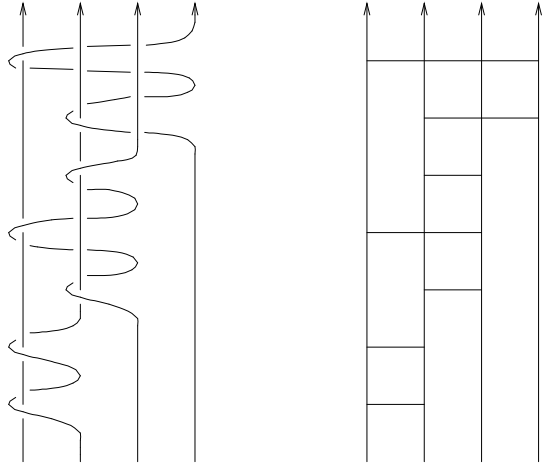

\[ \EEPIC{\CombedBraid}{0.6} \]
\caption{Left: In a combed braid, the first strand is always straight. The
second twists around the first some number of times, and then becomes
straight. The third waits patiently until the second finishes twisting
around the first, tangles a bit between the first and the second, and then
becomes straight. The fourth waits for the third to be done, tangles, etc.
Right: a non-decreasing chord diagram.}
\lbl{CombedBraid}
\end{figure}

\begin{pf}
\begin{enumerate}

\item Every braid can be `combed' by induction. Assume the first four
  strands are already combed as in figure~\ref{CombedBraid}, and that a
  fifth strand tangles between them. Think of the first four strands as
  made of copper wires and grease them very well. Think of the fifth
  strand as made of soft spaghetti, and use a fan to blow strong wind
  from the bottom up. The spaghetti wants to fly up, and the copper
  wires are very smooth (and sloped upwards) so they can't stop it from
  doing so. When all the spaghetti (except the very beginning, which is
  tied to the bottom plane) reaches the top of the figure, freeze it in
  place. It is now a path in the the plane minus four points (the four
  other strands). Present it as a product of generators, replace it by
  a copper wire, and you are ready to deal with the sixth strand. For a
  formal proof, see e.g.~\cite{Birman:BraidsBook}.

\item Immediate from fact~\ref{NonDecGen} (that non-decreasing chord
diagrams generate $\calA^{pb}$) and from the main technical observation of
this paper, proposition~\ref{MainTech}.

\end{enumerate}
\end{pf}

\begin{exercise} The (formal) Knizhnik-Zamolodchikov connection
\cite{Knizhnik-Zamolodchikov, Kohno:MonRep, Kontsevich:Vassiliev,
Bar-Natan:Vassiliev} defines a `holonomy map' $Z:K_n\to\calA^{pb}_n$.
Use simple properties of $Z$ (or of any of the other `universal
Vassiliev invariants' constructed in \cite{AltschulerFreidal:Universal,
Bar-Natan:NAT, Cartier:Construction, Kassel:QuantumGroups,
LeMurakami:Universal, Piunikhin:Combinatorial}) and
the obvious similarity between equation~\ref{SemiDirectF} and
equation~\ref{SemiDirectFA} to show that $Z$ is injective. Composing with
linear functionals on $\calA^{pb}_n$, one gets exactly the Vassiliev
invariants.
\end{exercise}
\section{Corollaries} \lbl{cors}

\subsection{The HOMFLY polynomial and braids}

It is well known (see
e.g.~\cite{Reshetikhin:Quasitriangle,ReshetikhinTuraev:Ribbon} that
$gl(N)$ in its defining representation corresponds to the HOMFLY
polynomial (in some parametrization), and that the higher
representations of $gl(N)$ correspond to various cabling operations
applied to the HOMFLY polynomial (for a similar situation,
see~\cite{MortonStrickland:Satellites}).

\begin{corollary} A pure braid is determined by the HOMFLY polynomials
of all closures into links of all of its cablings. (We allow an arbitrary
permutation of the strands before closing).
\end{corollary}

\begin{remark} \lbl{nonpurebraids}
The above corollary remains true even if the words `pure braid' are
replaced by the word `braid'. Indeed, it is easy to read the number of
components of a link from its HOMFLY polynomial, and knowing this
number for all possible closures of a braid (using all possible
permutations of the strands) determines the permutation $\sigma$
underlying that braid, as that number is maximal only if the closure is
done using the permutation $\sigma^{-1}$.
\end{remark}

\subsection{Braids and string links}

The Reshetikhin-Turaev $gl(N)$ invariants $J_{R_1,\ldots,R_n}$ were
originally defined as tangle invariants, and, in particular, they
extend to $n$-component string links (for a definition,
check~\cite{Habegger-Lin:Classification}
or~\cite{Bar-Natan:Homotopy}). Clearly, the same holds for the various
traces of the $J_{R_1,\ldots,R_n}$'s, and as these span the space of all
Vassiliev invariants of pure braids, it follows that

\begin{corollary} A Vassiliev invariant of pure braids can always be
extended to become a Vassiliev invariant (of the same type) of string
links.
\end{corollary}

On the level of chord diagrams (which are dual to invariants), this
corollary implies that chord diagrams of pure braids inject into chord
diagrams of string links. Recall that the latter space, $\calA^{sl}$,
discussed in more detail in~\cite{Bar-Natan:Homotopy}, is the space of
all diagrams of the form
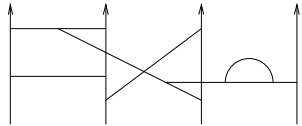
\[ \EEPIC{\Asl}{0.5} \]
(non-horizontal chords and oriented internal trivalent vertices are
allowed, apparent quadrivalent vertices in the planar projection are not
true vertices), modulo the STU relation,
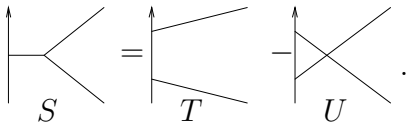
\[ \EEPIC{\STU}{0.5}. \]

\begin{corollary} \lbl{curiosity}
The obvious map $\calA^{pb}\rightarrow\calA^{sl}$ is an injection.
\end{corollary}

In fact, curiosity whether corollary~\ref{curiosity} holds is what lead me
to the investigations described in this paper.
\ifx\undefined\bysame
        \newcommand{\bysame}{\leavevmode\hbox to3em{\hrulefill}\,}
\fi

\end{document}